\documentclass{amsart}
\usepackage{times,amssymb,amsmath,amsfonts,amsthm}

\theoremstyle{plain}
\newtheorem{theorem}{Theorem}[section]
\newtheorem{proposition}[theorem]{Proposition}

\newtheorem{corollary}[theorem]{Corollary}
\theoremstyle{definition}
\newtheorem{example}{Example}[section]
\newtheorem{definition}{Definition}[section]
\theoremstyle{remark}
\newtheorem{remark}{Remark}[section]
\newtheorem*{acknowledgements}{Acknowledgements}

\newcommand{\id}{\operatorname{Id}}

\newcommand{\ob}{\operatorname{ob}}
\newcommand{\mor}{\operatorname{mor}}
\newcommand{\hot}{\operatorname{hot}}
\newcommand{\spec}{\operatorname{spec}}

\renewcommand{\ker}{\operatorname{Ker}}
\newcommand{\im}{\operatorname{Im}}
\newcommand{\coder}{\operatorname{Coder}}

\newcommand{\ds}{\displaystyle}

\newcommand{\Vect}{\mathfrak{Vect}}
\newcommand{\Mod}{\mathfrak{Mod}_\Lambda}
\newcommand{\Fun}{\mathfrak{Fun}}

\renewcommand{\L}{\mathfrak{L}}

\newcommand{\C}{\mathfrak{C}}

\newcommand{\Chat}{\widehat{\mathfrak{C}}}

\newcommand{\QB}{\mathfrak{B}}

\newcommand{\Sets}{\mathfrak{Sets}}

\newcommand{\F}{\mathbf{F}}
\newcommand{\T}{\mathbf{T}}
\newcommand{\D}{\mathbf{D}}
\newcommand{\G}{\mathbf{G}}
\newcommand{\un}{\mathbf{un}}
\newcommand{\fr}{\mathbf{fr}}
\newcommand{\z}{\mathbf{z}}

\newcommand{\Def}{\operatorname{\mathbf{Def}}}
\newcommand{\QFT}{\operatorname{\mathbf{QFT}}}
\newcommand{\MC}{\operatorname{\mathbf{MC}}}
\newcommand{\QM}{\operatorname{\mathbf{QM}}}
\newcommand{\CV}{\mathrm{C}(V)}

\newcommand{\Z}{\mathbb{Z}}

\newcommand{\m}{\mathfrak{m}}
\newcommand{\M}{\mathcal{M}}
\newcommand{\Y}{\mathcal{Y}}

\newcommand{\<}{\langle}
\renewcommand{\>}{\rangle}

\begin{document}

\title{Quantum backgrounds and $\QFT$}

\author{Jae-Suk Park}
\address{Dept. of Mathematics, Yonsei University,
    Seoul 120-749}
\thanks{The first author supported by KOSEF
    Interdisciplinary Research Grant No. R01-2006-000-10638-0.}
\email{jaesuk@yonsei.ac.kr}
\author{John Terilla}
\address{Queens College, City University of New York, 
    Flushing NY 11367}
\thanks{The second author supported in part by a grant from The
    City University of New York PSC Research Award.}
\email{jterilla@qc.cuny.edu}
\author{Thomas Tradler}
\address{NYC
    College of Technology, City University of New York, 
    Brooklyn NY 11201} 
\email{ttradler@citytech.cuny.edu}

\begin{abstract}We introduce the concept of a quantum background and 
a functor $\QFT$.  In the case that the $\QFT$ moduli
space is smooth formal, we construct 
a flat quantum superconnection on a bundle over $\QFT$
which defines algebraic structures relevant to
correlation functions in quantum field theory.  We
go further and identify 
chain level generalizations
of correlation functions which should be present in all quantum field
theories.  
\end{abstract}


\maketitle

\setcounter{tocdepth}{1}
\tableofcontents

\section{Introduction}\label{Introduction}
The work described in this paper grew from an effort to understand
quantum field theory mathematically.  In particular, we sought to
understand the way in which deformations of the action functional of a
quantum field theory control the correlation functions of the theory.
Our efforts led us to consider what we call \emph{quantum
  backgrounds}, or just \emph{backgrounds} for short.  Quantum
backgrounds were invented by the first author and introduced during a
sequence of lectures \cite{js2} in 2004 extending his work on flat
families of quantum field theories \cite{js1}.  We believe they
provide the proper algebraic setting in which local deformations of
quantum field theories should be discussed.  A background $B$ is a
four tuple $B=(P,m,N,\varphi)$ where $P$ is a graded noncommutative
ring, $m$ is a degree one element of $P$ satisfying $m^2=0$, $N$ is a
graded left $P$ module, and $\varphi$ is a degree zero element of $N$
satisfying $m\varphi=0.$ Additionally, we require that $P$ and $N$ be
free $k[[\hbar]]$ modules for a fixed field $k$, and that $P$ be a
$k[[\hbar]]$ algebra with $P/\hbar P$ commutative.  The condition
$m\varphi=0$ is analogous to having specified an initial action
functional that satisfies the quantum master equation.  In this paper,
we develop our ideas about backgrounds alongside a condensed overview
of deformation theory; we do this for two reasons.

The first reason is that physical information (such as correlation
functions, our generalized correlation functions, and generalized Ward
identities) are extracted from the background data with the aid of a
set valued functor that is comparable to the classical deformation
functor.  Essentially, we are interested in the equation
\begin{equation}\label{muphi}
  mU\varphi=0
\end{equation}
where $U$ is an invertible element of $P$ that depends on parameters.
We supply two interpretations: a Schr\"odinger interpretation where
$U\varphi$ is regarded as an evolution of $\varphi$, and a Heisenberg
interpretation where $U^{-1}mU$ is regarded as an evolution of $m$.
From either perspective, the evolution $U$ is not arbitrary and the
way it is constrained is intimately tied to the physics of the theory.
If the evolution is unobstructed, the constraints provides enough
relations to determine the correlation functions up to finite
ambiguity.  While one might think of $U\phi$ as a deformation of
$\phi$, or $U^{-1}mU$ as a deformation of $m$, we prefer the
interpretation as an evolution since that concept, along with the
definition of a quantum background itself, is evocative of a sort of
\emph{graded} quantum mechanics \`a la Dirac \cite{Dirac}, but it is a
quantum mechanics in the moduli space of action functionals for
quantum field theories having the same fields.

In this paper, we make restrictions on the form of $U$ and the type of
parameter space over which $\varphi$ evolves.  Recall the
Maurer-Cartan equation $d\gamma+\frac{1}{2}[\gamma, \gamma]=0$ in a
differential graded Lie algebra $L$ leads to a covariant set valued
functor of parameter rings $\Def_L$ where $\Def_L(A)$ is the set of
solutions to the Maurer-Cartan equation with parameters in the ring
$A$, modulo a natural equivalence.  In the same way, the quantum
master equation $mU\varphi=0$ for a background $B$ leads to a
set-value functor of $k[[\hbar]]$ parameter rings $\QFT_B$.  Standard
deformation theory, enriched to account for $\hbar$, is used to prove
(Theorem \ref{main_thm}) that the functor $\QFT_B$ is (pro, homotopy)
representable by a differential algebra $R$ and a versal solution $U$
to Equation \eqref{muphi} with parameters in $R$.  Recall that a
functor $\F$ is representable by a ring $R$ if it is naturally
equivalent to the functor $\hom(R,\,)$, in which case the identity
homomorphism $\id:R\to R$ corresponds to an element $U$, which is a
versal element in $\F(R)$, modulo the obstructions encoded by the
differential in $R$.  Here, the evolution of $\varphi\leadsto
U\varphi$ resembles a Schr\"odinger wave function over the space of
theories; or $U\varphi$ might be compared to a versal action
functional over the space of action functionals.

The main construction in the paper applies when $\QFT_B$ is smooth
formal, which in algebraic terms means that the differential in the
representing algebra $R$ vanishes; equivalently, there exist no
obstructions to constructing a versal action functional (see
Definition \ref{smoothformal}). In geometric terms, smooth formal
means the initial background data comprises a smooth point in the
formal moduli space.  In this case, we construct (Section
\ref{Flat_Connection}) a quantum flat superconnection on a bundle of over the
moduli space
$$\nabla:D_\Pi\to D_\Pi \otimes_R \Omega $$
Here, $D_\Pi$, $R$, and are the the tangent bundle, ring of functions,
and space of differential forms on the moduli space of solutions to
Equation \eqref{muphi}, all are free $k[[\hbar]]$ module equivalents
of the usual notions in formal geometry.  We call the map $\nabla$ a
quantum superconnection because: (i) $\nabla$ depends on $\hbar$ and
satisfies the $\hbar$-connection equation
\begin{equation}\label{conn_eq}
  \nabla(fe)=\hbar(df)e+ (-1)^{|f|} f\nabla(e) \text{ for $f\in R$ and
    a section $e$};
\end{equation}
and is a Quillen-type superconnection since the target involves all of
$\Omega$ and not only $\Omega^1$.  We prove (Theorem \ref{conn_thm})
that $\nabla^2=0$; i.e., $\nabla$ is flat.  The connection one-form
part of $\nabla$, call it $\nabla^{1}:D_\Pi \to D_\Pi \otimes_R
\Omega^1$, defines a flat, torsion free connection, which taken
together with $\hbar$ in Equation \eqref{conn_eq}, is similar to the
pencil of flat connections on the tangent bundle of a moduli space of
topological conformal field theories \cite{Dub}.  But $\nabla^{1}$
will, in general, depend on $\hbar$.  There is some evidence for the
existence of special coordinates \cite{js3,jt2} in which the $\hbar$
dependence of $\nabla^{1}$ reduces it to a pencil of flat connections,
equipping the moduli space with the structure of a representation of
the operad $H(\overline{\mathcal{M}}_{0,n})$, but this is not proved
in the present paper.  One can think of the flat quantum
superconnection $\nabla$ as furnishing a sort of Frobenius infinity
manifold structure to the moduli space associated to a background.

In the physics language, $\nabla^{1}$ encodes, up to a choice of one
point functions, all of the $n$ point correlation functions for $n>1$
of all the theories in the moduli space.  We have the compelling
metaphor: a theory of graded quantum mechanics over the moduli space
of quantum field theories is tantamount to the correlation functions
of the theories in the space.  This is accomplished without defining a
measure for the path integral.  From our point of view, the path
integral amounts to nothing more than a choice of the one point
functions, essentially a choice of basis of observables.  The
connection $\nabla^{1}$ determines the two point functions for every
theory in a neighborhood of the moduli space, and the $n$-point
functions for $n>2$ are obtained via products and derivatives in the
moduli directions.  The rest of the superconnection, $\nabla^k$ for
$k>1$, defines a homotopy algebraic structure which is a chain level
generalization of correlation functions which, to our knowledge, is
new in physics.  We offer an analogy with differential graded
algebras: compare the correlation functions derived from $\nabla^1$ to
the associative product in the homology of a differential graded
algebra and compare our chain-level generalized correlation functions
derived from $\nabla^k$ for $k>1$ to a minimal $A_\infty$ structure
defined in the homology of a dga that is quasi-isomorphic to the
original dga.  We suspect that the quantum background is determined up
to quasi-isomorphism by the flat quantum superconnection $\nabla$, but
we do not pursue this idea in this paper because $\nabla$ is defined
here only for smooth formal backgrounds and so it's difficult to
determine the extent to which $\nabla$ can be interpretted as a
``minimal model'' for the background until obstructed backgrounds are
better understood.  Also, we acknowledge that our definition of
``quasi-isomorphism'' for quantum backgrounds may need some
refinement---in particular, may need to encompass the concept of
special coordinates.

The second reason that we develop $\QFT$ alongside $\Def$ is our sense
that $\QFT$ generalizes $\Def$.  The relationship is epitomized in
Section \ref{Application}, where we discuss an example related to
differential BV algebras.  The deformation theory of quantum field
theories (i.e., the functor of $k[[\hbar]]$ algebras $\QFT_B$ for a
background $B$) is richer than classical deformation theory (i.e., the
functor of $k$-algebras $\Def_L$ for an $L_\infty$ algebra $L$).  The
fact that when $\QFT_B$ is a smooth functor, the moduli space admits
the quantum flat connection encoding correlation functions and
homotopy correlation functions is evidence that $\QFT_B$ is richer
since the $\Def_L$-moduli space does not naturally carry such an
additional structure.  The definition of a background is quite general
and might be constructed from a wide variety of mathematical data.
Our attitude is that the functor $\QFT_B$ provides a notion of quantum
deformation theory that in some examples extends the classical
deformation theory in the $\hbar$ direction.  When it is possible, one
should attribute a background $B$ to some initial mathematical data
under investigation, rather than an $L_\infty$ algebra $L$ defined
over $k$, and study functors of $k[[\hbar]]$-algebras.  Our
expectation is that correlations, at least in special coordinates,
arising from $\QFT_B$ reveal invariants of the initial data that is
more subtle than what can be extracted from $\Def_L$ as a functor of
$k$-algebras.  However, in order to make a precise statement about
$\QFT$ generalizing $\Def$, one should have a theory of special
coordinates and perhaps generalize even further the kind of parameter
rings on which $\QFT_B$ can be defined.  One might think of $\Def_L$
as the functor that controls deformations of objects over (formal,
graded, differential) commutative spaces whereas $\QFT_B$ controls
deformations of objects over something like ``quantized spaces.''  We
will expand on this point of view in another paper.

\begin{acknowledgements}
  We thank Michael Schlessinger, James Stasheff, and Dennis Sullivan
  for many helpful and inspiring conversations and, in particular to
  Dennis Sullivan, for generous and supportive working conditions at
  CUNY.  The authors also thank the referee for excellent
  suggestions: what clarity there is in the description of the
  flat quantum superconnection is due to the referee's recommendations,
  and while we were not able, in the time allotted, to convert
  ourselves to the
  ``functors from rings to 
  simplicial sets'' point of view, we recognize that it is a very good
  suggestion we intend to pursue it.
\end{acknowledgements}

\section{Categories, functors, and classical deformation theory}
\label{Primer}

This section is essentially a technical primer about deformation
theory approached via functors of parameter rings.  The reader may wish
to skip directly to section \ref{Quantum_Backgrounds} where our new
ideas are introduced.

\subsection{Categories}\label{Categories}

First, we review some of the tools we employ.  In order to describe
the moduli space of quantum backgrounds, we use functors of
parameter rings that are differential $k[[\hbar]]$-algebras, in
addition to being differential graded local Artin rings.

\subsubsection{Parameter rings}
We will be working with certain parameter rings, which we now
describe.  Fix a field $k$.  Although we have interest in the case
that $k$ has characteristic $p$, let us assume for this paper that the
characteristic of $k$ is zero.  A graded Artin local algebra is
defined to be a graded commutative associative unital algebra over $k$
with one maximal ideal satisfying the ascending chain condition.  Such
algebras are those local algebras whose maximal ideals are nilpotent
finite dimensional graded vector spaces.  We use $\m_A$ to denote the
maximal ideal of an Artin local algebra $A$.  The quotient $A/\m_A$ is
a field called the residue field of $A$.  A differential on $A$ is
defined to be a degree one derivation $d:A \to A$ satisfying $d^2=0$.
We denote the homology of $A$ by $H(A):=\ker(d)/\im(d)$.

Let $\C$ be the category whose objects are differential graded Artin
local algebras with residue field $k$ such that $[1]\neq [0]$ in
$H(A)$; morphisms are local homomorphisms, i.e. differential graded
algebra maps $A\to B$ such that $\im(\m_A)=\m_B$ and the induced maps
on residue fields is the identity.  Let $\Chat$ denote the category
whose objects consist of differential graded complete noetherian local
algebras $A$ for which $A/(\m_{A})^j \in \C$ for every $j$ and whose
morphisms are local.  Note objects in $\Chat$ are colimits of
objects in $\C$: $A=\underleftarrow{\lim}\, A/(\m_{A})^j $.

Let $\Lambda\in \ob(\Chat)$, concentrated in degree zero with zero
differential, and denote the maximal ideal $\m_{\Lambda}$ by $\mu$.
We define the category $\C_{\Lambda}$ to be the subcategory of $\C$
consisting of differential graded Artin local $\Lambda$-algebras;
morphisms in $\C_\Lambda$ are local differential $\Lambda$-algebra
homomorphisms.  Let $\Chat_\Lambda$ be the category of projective
limits of $\C_\Lambda$.  We will be concerned primarily with
$\Lambda=k[[\hbar]]$ or $\Lambda=k$ itself (and $\C_\Lambda=\C$).

\begin{definition}
  We call objects in $\C_\Lambda$ parameter rings.
\end{definition}

We have functors (``free'', ``zero action'', and ``underlying'')
\begin{equation}\label{fzu}
  \fr:\Chat \to
  \Chat_\Lambda,\quad \z:\Chat \to \Chat_\Lambda, \text{ and }\un:\Chat_\Lambda \to
  \Chat
\end{equation}
defined by $\fr(A)=A\hat{\otimes}_k \Lambda$, $\z(A)=A$ where the
action of $\Lambda$ on $A$ is given by $\lambda a=0$ for $\lambda \in
\mu$ and $a\in A$, and $\un(A)=A$, the underlying Artin algebra with
the $\Lambda$ action forgotten.  Note that $\un\circ \z=\id:\Chat \to
\Chat$.

\subsubsection{Tensor and fiber products of objects}

For $A,B\in \ob(\Chat_\Lambda)$, we define $A\otimes B:=
A\widehat{\otimes}_\Lambda B\in \ob(\Chat_\Lambda)$ with grading
determined by $\deg(a\otimes b)=\deg(a)+\deg(b)$ and differential $d$
defined as $d(a \otimes b)=d(a)\otimes b + (-1)^{\deg(a)}a \otimes
d(b)$, for homogeneous $a\in A$, $b\in B$.

Let $A,B,C \in \ob(\C_\Lambda)$ and $\alpha\in \hom(A,C)$ and $\beta
\in \hom(B,C)$.  We define the fibered product $A\times_C B \in
\ob(\C_\Lambda)$ by $A \times_C B=\{(a,b):\alpha(a)=\beta(b)\}$.  In
$\Chat_\Lambda$, fibered products do not exist since $A \times_C B$
may fail to be noetherian.

\subsubsection{Homotopy equivalence of morphisms and small extensions}

Let us now recall the homotopy model for the interval $I$.

\begin{definition}
  Define $I\ob(\Chat_\Lambda)$ by $I=k[[u,v]]$, $\deg(u)=0$,
  $\deg(v)=1$, and $d(p(u)+q(u)v)=p'(u)v.$
\end{definition}

One cannot evaluate an arbitrary element of $I$ at particular values
of $u$ and $v$, but for the (non-Artin) subalgebra $I'$ consisting of
polynomials in $u$ and $v$ and one does have evaluation maps
$ev_i:I'\to k$ and $ev_i:\fr(I')\to \Lambda$, where $ev_i$ is
evaluation at $(i,0)$ for $i=0,1.$

\begin{definition}
  Two morphisms $\tau_0,\tau_1:A \to B$ in $\Chat_\Lambda$ are called
  homotopy equivalent, or homotopic, if there exists an Artin
  subalgebra $B'\subset B \otimes \fr(I')$ and a morphism
  $\tau(u,v):A\to B'$ in $\Chat_\Lambda$ satisfying
  $\tau(0,0):=ev_0\tau(u,v)=\tau_0$ and
  $\tau(1,0):=ev_1\tau(u,v)=\tau_1.$ As suggested by the terminology,
  homotopy equivalence is an equivalence relation.  We denote the set
  of homotopy equivalence classes of morphisms $A \to B$ by
  $\hot(A,B).$
\end{definition}

\begin{definition}For $\mathfrak{A}=\C_\Lambda,\Chat_\Lambda$, we
  define the category $[\mathfrak{A}]$ to be the homotopy category of
  $\mathfrak{A}$.  That is, the category defined by
  $\ob([\mathfrak{A}])= \ob(\mathfrak{A})$ and $\mor(A,B)=\hot(A,B)$.
\end{definition}

Let $\rho:A \to B$ be a surjective morphism in $\C_\Lambda$.
\begin{definition}We say $\rho$ is a small extension if $\ker(\rho)$
  is a nonzero ideal $J$ such that $\m_A J=0$.  We say that $\rho$ is
  an acyclic small extension if, in addition, $H(J)=0$.
\end{definition}

\subsection{Functors of parameter rings}\label{Functors}
We now recall some basic definitions and properties of set valued
functors of parameter rings.  For details see \cite{S,M2}.

\subsubsection{Representability}
Any covariant functor $\F:\C_\Lambda \to \Sets$, the category of sets,
can be extended naturally to a functor $\F:\Chat_\Lambda \to \Sets$ by
$F(A)=\underleftarrow{\lim}\, F(A/\m_A^j).$ A functor
$\F:\Chat_\Lambda \to \Sets$ satisfying
$F(A)=\underleftarrow{\lim}\,F(A/\m_A^j)$ is called continuous.  We
consider continuous, covariant functors $\F:\Chat_\Lambda \to \Sets$
satisfying $\F(k)=F(\Lambda)=\{\text{one point}\}$. 
Such functors form a cateogry
$\Fun_\Lambda$ whose morphisms are natural transformations.  In the
case that $\Lambda=k$, we will denote $\Fun_\Lambda$ simply by $\Fun$.

Precomposition with $\fr,\z,$ and $\un$ defined in \eqref{fzu} gives
functors (which we denote by the same letters)
\begin{equation}\label{fzu2}
  \fr:\Fun_\Lambda \to \Fun, \quad \z:\Fun_\Lambda
  \to \Fun \text{ and } \un :\Fun \to \Fun_\Lambda.
\end{equation}

\begin{definition}We define a couple for $\F\in \Fun_\Lambda$ to be a
  pair $(A,\xi)$ where $A\in \Chat_\Lambda$ and $\xi \in \F(A).$
\end{definition}

\begin{definition}For $R\in \Chat_\Lambda$, we define
  $h_R:\Chat_\Lambda \to \Sets$ by $h_R(A)=\hom_\Lambda (R,A)$, and
  define $[h_R]:[\Chat_\Lambda]\to \Sets$ by
  $[h_R(A)]=\hot_\Lambda(R,A).$
\end{definition}

\begin{definition}
  We say that $\F\in \Fun_\Lambda$ is represented by $R$ if $\F$ is
  isomorphic to $h_R$ for some $R\in \Chat_\Lambda$.  We say that
  $\F\in \Fun_\Lambda$ is homotopy represented by $R$ if $\F$ induces
  a functor $[\Chat_\Lambda]\to \Sets$ that is represented by $[h_R]$.
\end{definition}

In the original terminology \cite{S}, Schlessinger called
representable functors and couples ``pro-representable'' and
``pro-couples''---we drop the prefix.  Homotopy representable functors
were treated by Manetti in \cite{M2} and Schlessinger-Stasheff in
\cite{SS}.

\begin{remark}If $\F$ is represented by $R$, then under the
  identification $\hom(R,R)\simeq \F(R)$, there is an element $\Pi\in
  \F(R)$ corresponding to $\id\in \hom(R,R).$ The element $\Pi$ is
  versal in the sense that for all $A\in \Chat_\Lambda$, and for all
  $\xi \in \F(A)$, there exists a unique morphism $\tau:R \to A$ so
  that $F(\tau)(\Pi)=\xi.$ One may say that $\F$ is represented by the
  couple $(R,\Pi)$.

  Likewise, if $\F$ is homotopy represented by $R$, then under the
  identification $\hot(R,R)\simeq \F(R)$, there is an element $\Pi\in
  \F(R)$ corresponding to $[\id]\in \hot(R,R)$.  The element $\Pi$ is
  versal in the sense that for all $A\in \Chat_\Lambda$, and for all
  $\xi \in \F(A)$, there exists a morphism $\tau:R \to A$ unique up to
  homotopy, so that $F(\tau)(\Pi)=\xi.$ One may say that $\F$ is
  homotopy represented by the couple $(R,\Pi)$.
\end{remark}

\subsubsection{Smoothness and quasi-smoothness}

\begin{definition}We say $\F\in \Fun_\Lambda$ is smooth if and only if
  the map $\F(\tau):\F(A)\to \F(B)$ is surjective for all small
  extensions $\tau:A \to B$.
\end{definition}
\begin{definition}
  We say $\F\in \Fun_\Lambda$ is quasi-smooth if and only if the map
  $\F(\tau):\F(A)\to \F(B)$ is surjective for all acyclic small
  extensions $\tau:A \to B$.
\end{definition}

\subsubsection{The Schlessinger condition}
Let $A,B,C \in \ob(\C_\Lambda)$ and $\alpha\in \hom(A,C)$, $\beta \in
\hom(B,C)$ and consider the fibered product $A\times_C B$.  For any
functor $\F \in \Fun_\Lambda$, the versal property of fibered products
(of sets) gives a map
\begin{equation}\label{the_eta_map}
  \eta:\F(A \times_C B) \to \F(A)\times_{\F(C)} \F(B)
\end{equation}

\begin{definition}
  We say that $\F \in \Fun_\Lambda$ satisfies the Schlessinger
  condition if the map $\eta$ in equation \eqref{the_eta_map} is a
  bijection.
\end{definition}

Readers familiar with Schlessinger's work may recognize what we call
the Schlessinger condition as a version of his condition (H4) (see
theorem 2.11 in \cite{S}).  This condition has also been called the
\emph{Mayer-Vietoris property} \cite{Maz} because of the
relationship to abstract homotopy theory \cite{B1}.

\subsubsection{Homotopy equivalence}\label{simdefined}
Let $\G \in \Fun_\Lambda$ and $A\in \ob(\C_\Lambda)$, we define a
relation on the set $\G(A)$ by setting $\xi_0 \sim \xi_1$ if and only
if there exists an Artin subalgebra $A'\subset A\otimes \fr(I')$ and a
$\xi(u,v) \in \G(A')$ with $\xi(0,0)=\xi_0$ and $\xi(1,0)=\xi_1.$
We're using the shorthand $\xi(i,0)$ for $\G(ev_i)(\xi(u,v))$ where
$ev_i:A \otimes \fr(I') \to A$ is the evaluation map determined by
$(u,v)=(i,0)$ for $i=0,1$.  We refer the reader to \cite{M2} for
details.
\begin{proposition}\label{equiv_thm}
  If $\G$ is quasi-smooth and satisfies the Schlessinger condition,
  then $\sim$ is an equivalence relation and the functor $\F\in
  \Fun_\Lambda$ defined by $\F(A)=\G(A)/{\sim}$ induces a functor on
  the category $[\Chat_\Lambda]$.
\end{proposition}

\begin{proof}
  The nontrivial, but standard, part is proving that the homotopy
  relation is transitive. (see \cite{M2} theorem 2.8).
\end{proof}

\subsubsection{Tangent space and tangent module}
For $i\in \Z$, we define the parameter ring $E_i$ by
$E_i:=k[\epsilon_i]/\epsilon_i^2$ where $\epsilon_i$ has degree $-i$
and $d(\epsilon_i)=0$.
\begin{definition}
  Then, for any $\F\in \Fun_\Lambda$, we define the tangent space to
  $\F$ by
$$T_\F=\oplus_i T_{\F}^i:=\oplus_i \F(E_i).$$
\end{definition}
If $\F$ commutes with certain products of rings with trivial
multiplicative structures then $T_\F$ naturally has the structure of a
graded vector space.  Specifically, if
\begin{equation}\label{vs_objects}
  \F(A \times_k B)\simeq \F(A)\times \F(B)
\end{equation}
for all rings $A$ and $B$ with $\m_A^2=\m_B^2=0$ and $d(A)=d(B)=0$
(the so-called vector space objects in $\C_\Lambda$) then $\F(A)$ will
be a vector space.  For example, if $\F$ satisfies the Schlessinger
condition, then $T_\F$ is a vector space.  Also, and importantly, if
$\F=\G/{\sim}$ as in proposition \ref{equiv_thm}, then $T_\F$ is a
vector space.  Furthermore, $\T$ defines a functor from the
subcategory of $\Fun_\Lambda$ consisting of those functors satisfying
\eqref{vs_objects} to $\Vect$, the catgory of graded vector spaces.

One can augment the tangent space slightly and obtain a natural
$\Lambda$-module structure.  For each $\alpha \in \Lambda$, one uses
the morphisms
$$\oplus_i \fr(E_i)=\oplus_i \Lambda[\epsilon_i]/\epsilon_i^2 \overset{\alpha
  \cdot}{\longrightarrow} \oplus_i
\Lambda[\epsilon_i]\epsilon_i^2=\oplus_i \fr(E_i)$$ to get maps
$\alpha:\oplus_i \F(\fr(E_i)) \to \oplus_i \F(\fr(E_i))$.

\begin{definition}
  For any $\F\in \Fun_\Lambda$, we define the tangent $\Lambda$ module
  to $\F$ by
$$D_\F:=\oplus_i D_\F^i \oplus_i \F(\fr(E_i)).$$
\end{definition}

Like the tangent space, $\D$ defines a functor from the subcategory of
$\Fun_\Lambda$ consisting of those functors satisfying
\eqref{vs_objects} to $\Mod$, the catgory of $\Lambda$-modules.

\begin{proposition}\label{rep_thm}
  Suppose that $\G\in \Fun_\Lambda$ is quasi-smooth and satisfies the
  Schlessinger condition.  Let $\F=\G/{\sim}$ as in proposition
  \ref{equiv_thm} and suppose that $T_\F$ is finite dimensional with
  basis $\{\xi_i\}$.  Let $\{t_i\}$ be the dual basis for $T_\F$.
  Then, there exists a differential $\delta$ on the algebra
  $R:=\Lambda[[t_i]]$ with $\delta(\m_R)\subseteq (\m_R)^2$ and an
  element $\Xi\in \F(R,\delta)$ so that $\F$ is homotopy represented
  by the couple $((R,\delta),\Xi)$.
\end{proposition}

\begin{proof}
  See theorem 4.5 and corollary 4.6 in \cite{M2}.
\end{proof}

Let us recall one more result from the general theory of functors of
parameter rings, which provides a theoretical underpinning for the
definition of ``quasi-isomorphism'' in sections \ref{Def_Theory} and
\ref{Quantum_Backgrounds}.

\begin{proposition}\label{inv_fun_thm}
  Suppose $\G,\G' \in \Fun_\Lambda$ are quasi-smooth and satisfy the
  Schlessinger condition.  Let $F=\G/{\sim}$ and $F'=\G'/{\sim}$ and
  suppose $\mathbf{n}:\F\to \F'$ is a natural transformation.  Then
  $\mathbf{n}$ is an isomorphism of functors if and only if
  $D_{\mathbf{n}}:D_\F \to D_{\F'}$ is an isomorphism of graded
  $\Lambda$ modules.
\end{proposition}

\subsection{Classical deformation theory}\label{Def_Theory}
We now recall the basic elements of classical deformation theory.  In
this subsection we define the category of $L_\infty$ algebras and
define a functor $\Def_L$ for each $L_\infty$ algebra $L$.  At the end
of the section, we highlight the properties of $\Def_L$ (as described
in the section \ref{Functors} for general functors) as they pertain to
$L$.

Now, we set $\Lambda=k$.

\subsubsection{$L_\infty$ algebras}
Let $V=\oplus_{j\in \mathbb{Z}} V^j$ be a graded vector space over
$k$.  As usual, let $V[n]$ denote the shift
$V[n]=\oplus_{j\in\mathbb{Z}} V[n]^j$, where $V[n]^j=V^{j+n}$.  Let
$S^i V$ denote the $i$-th symmetric product of $V$.  Note that
$\wedge: S^i V[1] \otimes S^j V[1] \to S^{i+j} V[1]$ makes
$$\CV:=\oplus_{i=0}^\infty
\mathrm{S}^i(V[1])$$ into a differential graded commutative
associative algebra.  It is also a graded cocommutative coassociative
coalgebra.  The coproduct $\CV\to \CV\otimes \CV$ is characterized by
the property that it is map of differential graded algebras, and that
$x\mapsto x\otimes 1 + 1\otimes x$ for each $x\in V[1]$.  Indeed,
$C(V)$ is a construction of the free cocommuative coassociative
coalgebra over $V[1]$.  The versal property characterizing this free
construction is: for any graded coalgebra $C$ and any graded linear
map $C\to V[1]$, there exists a unique coalgebra map $C \to C(V)$
extending the linear map.  This gives an isomorphism of graded vector
spaces: $\hom_\mathrm{linear}(C,V)\simeq
\hom_{\mathrm{coalgebra}}(C,\CV).$

We mention a second identification.  Any homomorphism $\sigma :
S^iV[1] \to V$ extends uniquely to a map $\overline{\sigma}\in
\coder(\CV)$, and any coderivation $\sigma \in \coder(\CV)$ is
determined by its components $\underline{\sigma}_i :S^iV[1] \to V$.
Thus, $\hom_{\mathrm{linear}}(\CV,V)\simeq \coder(\CV).$

\begin{definition}
  An $L_\infty$ algebra $L$ is defined to be a pair $L=(V,D)$ where
  $V$ is a graded vector space, and $D\in \coder^1(\CV)$ is a
  coderivation of degree $1$, satisfying $D^2=0$.  A morphism
  $\sigma:L\to L'$ between two $L_\infty$ algebras $L=(V,D)$ and
  $L'=(V',D')$ consists of a degree zero coalgebra map $\sigma:\CV \to
  \mathrm{C}(V')$ with $\sigma\circ D=D'\circ \sigma$. Let us denote
  the category of $L_\infty$ algebras by $\L$.
\end{definition}

We will use the notation $d_i$ for $\underline{D}_i:S^iV[1]\to V$, the
$i$-th component of $D$.

\subsubsection{Maurer-Cartan Functor}

\begin{definition}Let $L=(V,D)\in \ob(\L)$. We define a functor
  $\MC_L\in \Fun$ by
  $$\MC_L(A)=\{\text{degree $0$ differential coalgebra
    maps }\m_A^* \to \CV\}.$$ Here, $\m_A^*$ is the differential
  coalgebra $\m_A^*:=\hom(\m_A,k)$.
\end{definition}

The functor $\MC_L$ satisfies the Schlessinger condition and is
quasi-smooth.  Therefore, we can form the quotient $\MC_L/{\sim}$,
where $\sim$ is the natural gauge equivalence defined in Section
\ref{simdefined}.  Let us make a list of definitions:

\begin{remark}\label{MC_remark}
  Any coalgebra map $\m_A^* \to \CV$ is determined by a linear map
  $\m_A^* \to V$; equivalently, by and element in $V\otimes \m_A$.
  Then
  \begin{equation}\label{MC_remark_eq}
    \gamma \in \MC_L \Leftrightarrow d_1(\gamma)+\frac{1}{2!}d_2(\gamma\wedge
    \gamma)+\frac{1}{3!}d_3(\gamma\wedge \gamma\wedge \gamma)+\cdots=0.
  \end{equation}
  where the $\gamma$ on the right is viewed as an element in $V\otimes
  \m_A$ and $d_i:S^i((V \otimes \m_a) [1] \to V\otimes \m_A$ is
  extended from the components of the $L_\infty$ structure on $\CV$
  using the differential and multiplication in $A$.
\end{remark}

\begin{definition} Let $L$ be an $L_\infty$ algebra.
  \begin{enumerate}
  \item Define the quotient $\Def_L:=\MC_L/{\sim}$ where $\sim$ is the
    natural gauge equivalence.
  \item The assignment $L\mapsto \Def_L$ defines a functor from $\Def:
    \L \to \Fun$.
  \item We define $H(L)$ to be the vector space $H(L):=T_{\Def_L}[1]$
    and call it the homology of $L$.  The composition $H:=\T \circ
    \Def$ defines a functor $H:\L \to \Vect$.
  \end{enumerate}
\end{definition}

\begin{remark}
  In practice, the classical deformation theory of a mathematical
  object \cite{S,K1,M2} is controlled by an $L_\infty$ algebra $L$.
  This means that given a mathematical object, for which there exists
  the concept of a deformation of that object over a parameter ring
  $A$ and the concept of equivalence of such deformations, there is a
  bijection of sets
$$
\left\{
  \begin{array}{c}
    \text{(equivalence classes of) deformations}\\\text{of the object over
      the ring $A$}
  \end{array}\right\}
\leftrightarrow \Def_L(A).
$$
\end{remark}

\begin{proposition}\label{universal_mc_thm}
  Suppose $L$ is an $L_\infty$ algebra and $T_{\Def_L}=H(L)[-1]$ is
  finite dimensional with basis $\{\alpha^1, \ldots, \alpha^k\}$.  Let
  $\{s_i\}$ be the dual basis and let $S:=k[[s_1,\ldots, s_k]].$ There
  exists a differential $d:S\to S$ with $d(\m_S)\subset (\m_S)^2$ and
  $\Gamma \in \MC_L(S)$ such that $\Def_L$ is homotopy represented by
  the couple $((S,d),\Gamma)$.
\end{proposition}

\begin{proof}Apply proposition \ref{rep_thm} to $\Def_L$.
\end{proof}

\begin{definition}Let $L$ be an $L_\infty$ algebra and assume $H(L)$
  is finite dimensional.
  \begin{enumerate}
  \item We call $L$ finite.
  \item A morphism $L \to L'$ is called a quasi-isomorphism if the map
    $\Def_{L} \to \Def_{L'}$ is an isomorphism of functors, in which
    case $L$ and $L'$ are called quasi-isomorphic.
  \item We call $\Gamma$ in proposition \ref{universal_mc_thm} a
    versal solution to the Maurer-Cartan equation.
  \item We call $L$ smooth formal if $d=0$ for $d:S\to S$ as in
    proposition \ref{universal_mc_thm}.
  \end{enumerate}
\end{definition}

\begin{remark}
  In light of the proposition \ref{inv_fun_thm}, one can assemble the
  parts of the definition above to say that $L$ and $L'$ are
  quasi-isomorphic if and only if there exists a morphism $\sigma:L
  \to L'$ inducing an isomorphism $H(L) \to H(L')$, which together
  with a computation establishing
$$T_{\Def_L}\simeq \ker({d_1})/\im({d_1}),$$
may be used to give an alternative definition of quasi-isomorphism as
a morphism $L\to L'$ inducing an isomorphism in homology.
\end{remark}

\section{Quantum backgrounds}\label{Quantum_Backgrounds}
Quantum backgrounds were invented to provide an algebro-mathematical
formulation of quantum field theory.  We begin by defining a quantum
background (see Definition \ref{background_def}) which one can think
of as the input data, leading naturally to a functor $\QFT$ (see
Definition \ref{QFT_definition}).  From now on, we set
$\Lambda=k[[\hbar]]$.

\subsection{The category $\QB$}

\begin{definition}
  We say that a graded associative, unital $\Lambda$-algebra $P$ has a
  classical limit provided $P$ is free as a $\Lambda$ module and if
  $K=P/\hbar P$ is a graded, commutative, associative, unital
  $k$-algebra.
\end{definition}

Note that if $P$ has a classical limit, then $P\simeq K\oplus \hbar K
\oplus \hbar^2 K \oplus \cdots$ and $[P,P]\subseteq \hbar P.$

\begin{definition}
  \label{background_def}
  We define a background $B$ to be a four-tuple $B=\left( P,N, m,
    \varphi \right )$ where
  \begin{enumerate}
  \item $P=\oplus_{i} P^i$ is a graded, associative, unital
    $\Lambda$-algebra with a classical limit,
  \item $N=\oplus_{i}N^i$ is a graded left $P$ module, which is free
    as a $\Lambda$ module,
  \item $m\in P^1$ satisfies $m^2=0$ ($m$ is called a structure),
  \item $\varphi \in N^0$ satisfies $m \cdot \varphi=0$ ($\varphi$ is
    called a vacuum).
  \end{enumerate}

  A morphism between two backgrounds $B=\left( P,m,N, \varphi \right
  )$ and $B'=\left( P',m', N', \varphi' \right )$ consists of a map
  $\sigma: P \to P'$ of graded $\Lambda$ algebras and a map $\tau:N\to
  N'$ of graded $P$ modules (where $N'$ becomes a $P$ module via
  $\sigma$) with $\tau(\varphi)=\varphi'$ satisfying the compatability
$$\tau(ma\varphi)=m'\sigma(a) \varphi' \text{ for all }a\in P.$$
Denote the category of backgrounds by $\QB$.
\end{definition}

\subsection{Quantum master equation}\label{QME-section}

\begin{definition}
  Let $B=(P,N, m, \varphi)$ be a background, let $(A,d) \in
  \ob(\C_\Lambda)$, where $\C_\Lambda$ is the category of differential graded Artin local $\Lambda=k[[\hbar]]$ algebras, and $\pi \in (P\otimes\m_A)^0$.  We call the equation
  \begin{equation}\label{QME}
    \left( e^{{\pi}/{\hbar}}me^{-{\pi}/{\hbar}}-\hbar
      e^{{\pi}/{\hbar}}d \left( e^{-{\pi}/{\hbar}}\right)\right)
    \varphi=0
  \end{equation}
  the quantum master equation, and denote the set of solutions by
  $\QM_B(A)$:
  \begin{equation*}
    \QM_B(A)=\left\{ \pi \in  (P\otimes \m_A)^0
      \text{ such that $\left(
          e^{{\pi}/{\hbar}}me^{-{\pi}/{\hbar}}-\hbar
          e^{{\pi}/{\hbar}}d \left( e^{-{\pi}/{\hbar}}\right)
        \right)   \varphi=0$}\right\}.
  \end{equation*}
  Notice that $\m_k=0$, so $\QM_B(k)=\{0\}$, and any morphism $\tau:A
  \to B\in \hom(\C_\Lambda)$ gives rise to a map from $\QM_B(A) \to
  \QM_B(B)$ given by $\pi \mapsto (1\otimes \tau) \pi.$ Thus, we have
  defined a functor $\QM_B\in \Fun_\Lambda.$
\end{definition}

Definition \ref{QME} relies on the fact that both terms
$e^{{\Pi}/{\hbar}} m e^{- {\Pi}/{ \hbar}}$ and $\hbar
e^{{\Pi}/{\hbar}} d\left( e^{-{\Pi}/{ \hbar}}\right)$ are
well defined elements of $P\otimes A$, a fact that we clarify in two
remarks.  An example follows in a third remark.

\begin{remark}
  First, we establish the meaning of $\frac{[\eta,\zeta]}{\hbar}$ in
  $P\otimes A$.  The ring $P$ by assumption is free as a $\Lambda$
  module, but $A$, owing to its nilpotency, cannot be and so $\hbar$
  is a zero divisor in $P\otimes A$.  However, in $P$, one has
  $[P,P]\subseteq \hbar P$, so one has an operator
  \begin{align*}
    \frac{[\,,\,]}{\hbar}:P \times P & \to P\\
    \alpha,\beta \mapsto \gamma
  \end{align*}
  where $\gamma \in P$ is defined by expressing $[\alpha,\beta]=\hbar
  \gamma$.  Then, for $\eta=\alpha \otimes a$ and $\zeta=\beta\otimes
  b$ in $P\otimes A$, one defines
  $$\frac{[\eta,\zeta]}{\hbar}=\frac{[\alpha \otimes a,\beta\otimes
    b]}{\hbar}:=\gamma \otimes ab\in P\otimes A.$$ Then,
  $\frac{1}{\hbar^j}[\eta_1,[\cdots [\eta_{j-1},\eta_j
  ]\cdots]\subseteq P\otimes A$ is understood in the same way.
\end{remark}

\begin{remark}\label{UdU-1}
  Expanding $e^{{\pi}/{\hbar}} m e^{-{\pi}/{ \hbar}}$ and
  $\hbar e^{{\pi}/{\hbar}} d\left( e^{- {\pi}/{
        \hbar}}\right)$ in terms of repeated commutators:
$$e^{{\pi}/{\hbar}} m e^{-{\pi}/{ \hbar}}=\exp\left({\frac{[-,\pi]}
    {\hbar}}\right) (m)=m+{\frac{[m,\pi]} {\hbar }}+ {\frac{[[m,\pi],
    \pi]} {2 \hbar^2}}+\cdots$$ displays each term as an element of
$P\otimes A$, and the sum terminates since $\m_A$ is nilpotent.  The
same holds for the $d$ term once expanded as
$$\hbar e^{{\pi}/{\hbar}} d \left( e^{- {\pi}/{ \hbar}}
\right)=f \left({\frac{[-,\pi]} {\hbar}}\right) (d \pi), \text{ where
}f(x)= \frac{e^x -1}{x}=\sum_{j\geq 0}\frac{x^j}{(j+1)!}.$$
\end{remark}

\begin{remark}
In certain cases, the quantum background $B$ arises from a differential
BV algebra $(V,d,\Delta,\cdot)$ with Lie bracket $(\,,\,).$ In these
cases, the quantum master
equation defined for the background as defined above 
reduces to the usual quantum master equation
seen in the physics literature:
$$d\pi+\hbar \Delta \pi+\frac{1}{2}(\pi,\pi)=0$$
for $\pi\in V\otimes \m_A$. 
This example is illustrated in Section
\ref{Application}.
\end{remark}

\begin{theorem}\label{quasismooth}
  The functor $\QM_B$ is quasismooth and satisfies the Schlessinger
  condition.
\end{theorem}

\begin{proof}
  Since $\QM_B$ is given as a solution set of an equation, $\QM_B$
  satisfies the Schlessinger condition.  To show that it is
  quasi-smooth, let $(A,d) \in \ob(\C_\Lambda)$ and let $C_A=(P\otimes
  A) \varphi$ denote the cyclic submodule of $N\otimes A$ generated by
  the vacuum $\varphi$ ($=\varphi\otimes 1)$.  Note that $C_A$ is a
  complex with differential $D:C_A^j \to C_A^{j+1}$ defined by
  $$D((\pi \otimes a) \varphi)= (m\pi\otimes a)\varphi+(-1)^{|a|}\hbar
  (\pi \otimes d(a)) \varphi.$$

  Let $\tau: A \to A'$ be an acyclic small extension with kernel $I$
  and let $a' \in \QM_B(A')$.  Since $\tau$ is surjective, we can
  choose an $a \in P\otimes \m_A$ so that $1 \otimes \tau (a)=a'.$ Let
  $$i=e^{{a}/{\hbar}}me^{-{a}/{\hbar}}-\hbar
  e^{{a}/{\hbar}}\delta \left( e^{-{a}/{\hbar}}\right)\in
  P_A.$$ Notice that $(1\otimes \tau)(i\varphi)=0$.  So, $i\varphi$ is
  in the kernel of $1\otimes\tau:C_A \to C_{A'}.$ Since the functor
  $\otimes$ is right-exact, $i\varphi \in C_{I}$, and we can write
  $i\varphi=i'\varphi$ for some $i'\in I.$ Now, since $a\in \m_A$, we
  have $a i'=0$ and $e^{-{a}/{\hbar}}i'=i'.$ Now, we have
  \begin{multline*}
    i'\varphi=e^{-{a}/{\hbar}}i'\varphi
    =e^{-{a}/{\hbar}}i\varphi\\
    =e^{-{a}/{\hbar}}\left(e^{{a}/{\hbar}}me^{-{a}/{\hbar}}-\hbar
      e^{{a}/{\hbar}}d \left( e^{-{a}/{\hbar}}\right)\right)
    \varphi =\left(me^{-{a}/{\hbar}}-\hbar d\left(
        e^{-{a}/{\hbar}}\right)\right) \varphi.
  \end{multline*}
  Now, consider $$D(i'\varphi)=\left(m^2e^{-{a}/{\hbar}}-\hbar md
    \left(e^{-{a}/{\hbar}}\right) + d m e^{-{a}/{\hbar}}-\hbar
    d^2 \left( e^{-{a}/{\hbar}}\right)\right) \varphi=0.$$ Since
  $I$ is acyclic, $C_{I}$ is acyclic, which implies that there exists
  a $j\varphi \in C_{I}^{0}$ with $Dj\varphi =i'\varphi.$ Define
  $a'=a-\hbar j\in (P\otimes \m_A)^0.$ The claim is that $a' \in
  \QM_B(A).$ To show this, we compute.  First note that $aj=0$ and
  $j^2=0$ imply that $e^{{a'}/{\hbar}}=e^{{a}/{\hbar}}+j$ and
  $e^{{-a'}/{\hbar}}=e^{{-a}/{\hbar}}-j$.
  So
  \begin{multline*}
    \biggr(e^{{a'}/{\hbar}}me^{-{a'}/{\hbar}}-\hbar
    e^{{a'}/{\hbar}} d \left(
      e^{-{a'}/{\hbar}}\right)\biggr) \varphi\\
    =
    \left(\left(e^{{a}/{\hbar}}-j\right)m\left(e^{-{a}/{\hbar}}+j\right)-\hbar
      \left(e^{{a}/{\hbar}}-j\right)d \left(
        e^{-{a}/{\hbar}}+j\right)\right)\varphi \\
    =i\varphi-jm\varphi-mj\varphi-jmj\varphi+\hbar j d
    e^{-{a}/{\hbar}}\varphi -\hbar e^{{a}/{\hbar}} d j \varphi
    =(i'-Dj )\varphi =0.
  \end{multline*}

\end{proof}

As a corollary, we make the following definition:

\begin{definition}\label{QFT_definition}Let $B$ be a background.
  \begin{enumerate}
  \item We define $\QFT_B:=\QM_B/{\sim}$ where $\sim$ is the natural
    gauge equivalence.
  \item The assignment $B\mapsto \QFT_B$ defines a functor from $\QFT:
    \QB \to \Fun_\Lambda$.
  \item We define $H(B)$ to be the vector space $H(B):=T_{\QFT_B}$ and
    call it the homology of $B$.  The composition $H:=\T \circ \QFT$
    defines a functor $H:\QB \to \Vect$.  In addition, we define the
    Dirac module of $B$ to be the $\Lambda$-module $D(B):=D_{\QFT_B}.$
    The composition $D:=\D \circ \QFT$ defines a functor $D:\QB \to
    \Mod$.
  \end{enumerate}
\end{definition}

\begin{proposition}\label{H(B)-and-D(B)}
  Let $B$ be a background.  
There are natural isomorphisms
  \begin{gather}
    D(B)\simeq \left \{a\in P :
      \frac{[m,a]}{\hbar}\varphi=0\right\}\biggr/ \sim
    \label{dirac_computation} \intertext{and} H(B)\simeq \left \{a\in
      P \otimes_\Lambda k :
      \frac{[m,a]}{\hbar}\varphi=0\right\}\biggr/ \sim
    \label{tangent_computation}
  \end{gather}
  where $a_1 \sim a_0$ if and only if $(a_1-a_0)\varphi = \frac{[m,b]}{\hbar}
  \varphi$ for some $b\in P$ (or $P\otimes_\Lambda k$).
\end{proposition}
\begin{proof}
  We set up the natural isomorphism \eqref{tangent_computation} in
  homogeneous pieces by
$$\xi \leftrightarrow a.$$
For $\xi \in (P \otimes_\Lambda m_{E_i})^0$, with $\xi=a \epsilon_i$,
where $a\in P^{-i}\otimes_\Lambda k $ and $\epsilon_i \in m_{E_i}$,
one writes
$$e^{{a\epsilon_i}/{\hbar}}me^{-{a\epsilon_i}/{\hbar}} \varphi = \left(
  m-\frac{1}{\hbar}[m,a]\epsilon_i+\frac{1}{2
    \hbar^2}[[m,a],a]\epsilon_i^2+\cdots\right)\varphi =
-\frac{1}{\hbar}[m,a]\epsilon_i\varphi$$ and sees that $\xi \in
\QM_B(E_i)$ if and only if $\frac{[m,a]}{\hbar}\varphi= 0$.

Now we check equivalence.  First suppose
$\frac{[m,a_0]}{\hbar}\varphi= 0$ and $\frac{[m,a_1]}{\hbar}\varphi=
0$ and $(a_1-a_0)\varphi = \frac{[m,b]}{\hbar}\varphi$ for some $b\in
P.$ Then, define $\xi(u,v)=(a_1u+a_0(1-u)+bv)\epsilon_i.$ One sees
that $\xi(0,0)=\xi_0=a_0\epsilon_i$ and
$\xi(1,0)=\xi_1=a_1\epsilon_i$.  To see that $\xi(u,v)\in
\QM_B(E_i\otimes \fr(k[u,v]/(u^2,uv,v^2))$, compute:
\begin{multline*}
  \left(e^{{\xi(u,v)}/{\hbar}}me^{-{\xi(u,v)}/{\hbar}} -\hbar
    e^{{\xi(u,v)}/ {\hbar}}\delta\left(e^{-{\xi(u,v)}/{\hbar}}\right)
  \right)\varphi \\ =- \left( \frac{[m,\xi_1]}{\hbar}u
    \frac{[m,\xi_0]}{\hbar}(1-u)+\frac{[m,b]}{\hbar}v-(\xi_1-\xi_0)v\right)\varphi
  =0.
\end{multline*}
Conversely, let
$$\xi(u,v)=(a(u)+b(u)v)\epsilon_i \in P^{-i}\otimes I'$$
for $a(u), b(u)\in P^{-i}\otimes_\Lambda k[u]$ for some $I' \subseteq
k[u,v]$.  Then $\xi(u,v)\in \QM_B(E_i\otimes \fr(I'))$ implies that
$\frac{m,[a(u)]}{\hbar}\varphi=0$ and
$a'(u)\varphi=\frac{[m,b(u)]}{\hbar}\varphi$.  Define $b=\int_0^1
b(u)du$.  Then, for $\xi(0,0)=a_0\epsilon_i$ and
$\xi(1,0)=a_1\epsilon_i$, $a_0,a_1\in P^{-i}\otimes_\Lambda k$, we
have $\frac{[m,a_0]}{\hbar}\varphi= 0$ and
$\frac{[m,a_1]}{\hbar}\varphi= 0$ and $(a_1-a_0)\varphi =
\frac{[m,b]}{\hbar}\varphi$.

The computation for $D(B)$ is the same, except that the correspondence
\eqref{dirac_computation} $[\xi] \leftrightarrow [a]$ is setup for
$\xi\in (P\otimes_\Lambda m_{\fr(E_i)})^0 =P^{-i} \otimes_k m_{E_i}$
with $\xi=ae_i$ with $a\in P^{-i}$ and $e_i\in m_{E_i}$.
\end{proof}

\begin{corollary}
Let $C_\Lambda=P\varphi$ and $C_k=(P\otimes_\Lambda k)\varphi$ be the
  cyclic submodules of $N$ and $N\otimes k$ 
  generated by the vacuum.
  Then multiplication on the left by $m$ defines differentials 
  $m:C_\Lambda ^j \to C^{j+1}$ and $m:C_k^j\to C_k^{j+1}$.  Then the
  Dirac space and tangent space are isomorphic to the homology of
  these complexes
 $D(B)\simeq H(C_\Lambda ,m)$ and $H(B)\simeq H(C_k,m)$.
\end{corollary}

\begin{theorem}\label{main_thm}
  Suppose that $B$ is a background and $T_{\QFT_B}=H(B)$ is finite
  dimensional.  with basis $\{\xi^1, \ldots, \xi^r\}$.  Let
  $\{t_1,\ldots, t_r\}$ be the dual basis and let
  $R=\Lambda[[t_1,\ldots, t_r]]$.  There exists a differential
  $\delta:R \to R$ with $\delta(\m_R)\subset (\m_R)^2$ and $\Pi \in
  \QM_B(R)$ such that $\QFT_B$ is homotopy represented by the couple
  $((R,\delta),\Pi).$
\end{theorem}

\begin{proof}Apply proposition \ref{rep_thm} to $\QFT_B$.
\end{proof}

\begin{definition}\label{smoothformal}
  Let $B$ be a background and assume $H(B)$ is finite dimensional.
  \begin{enumerate}
  \item We call $B$ finite.
  \item A morphism of backgrounds $B \to B'$ is called a
    quasi-isomorphism if $\QFT_{B} \to \QFT_{B'}$ is an isomorphism of
    functors, in which case we say $B$ and $B'$ are quasi-isomorphic.
  \item We call $\Pi$ in proposition \ref{main_thm} a versal solution
    to the quantum master equation.
  \item We call $B$ smooth formal if $\delta=0$ for $\delta$ in
    proposition \ref{main_thm}.
  \end{enumerate}
\end{definition}

\section{Flat quantum superconnection} \label{Flat_Connection} One is
able to supply derived algebro-geometric interpretations of structures
associated to smooth formal backgrounds.  Some of these structures
have been unearthed in specific examples, and our aim in this section
to explain this structure from a unifying perspective afforded by
backgrounds.  For example, our main construction is a flat
superconnection $\nabla$ on a bundle over the functor $\QFT_B$.  The
connection one-form part of $\nabla$, which we call
$\nabla^{1}$, defines a flat connection which coincides with the
connection defined in \cite{B2} in the setting of quantum periods.
Throughout this section we assume that $B$ is a smooth formal background.

\subsection{Conceptual description}

First, to define a bundle over a functor $\F$ of parameter rings, one
assigns an $A$ module $M_\pi$ to each point $\pi\in \F(A)$, functorial
with respect to ring maps $A\to A'$.  In the case of $\QM_B$ for a
background $B$ and a point $\pi\in \QM_B(A)$ where $A$ is a parameter
ring with zero differential, we can define a new background $B_\pi$ as
$$B_\pi:=(B\otimes A, m^\pi, N\otimes A, \varphi \otimes 1)$$
where $m^\pi:=e^{\pi/ \hbar}me^{-\pi/ \hbar}$.  To see that $B_\pi$
defines a background, note that $m^\pi$ squares to zero (trivially)
and that $m^\pi (\varphi \otimes 1)=0$ by the quantum master equation.
Let
\begin{align*}
  D_\pi&:=D(B_\pi) = \{x \in P\otimes A: \frac{\left[m^\pi,
      x\right]}{\hbar} \phi =0\}/\sim \\
  &=H(C_A,m^\pi)
\end{align*}
where $C_A=(P\otimes A)\varphi$ is the cyclic submodule of $N\otimes
A$ generated by the vacuum.
The Dirac space $D_\pi$ of the background $B_\pi$ is an $A$ module.
If $f:A\to A'$ is a map of parameter rings, we obtain a map of
background $B_\pi \to B_{\F(f)(\pi)}$ and hence a map of Dirac spaces.
Furthermore, if $\pi\sim\pi'$, the backgrounds $B_\pi$ and $B_{\pi'}$
are quasi-isomorphic, hence have isomorphic Dirac spaces.  Thus, the
construction of $D_\pi$ from $\pi \in \QM_B(A)$ defines a bundle over
the functor $\QFT_B$ when $\QFT_B$ is smooth formal.  Call this bundle
the versal Dirac bundle and denote it by $\widetilde{D}$.  We haven't
defined a background $B_\pi$ for $\pi \in QM_B(A)$ when $A$ is a
parameter ring with nonzero differential, but it's straightforward
(but not necessary for smooth backgrounds) to modify the definition of
$D_\pi$ for such a point $\pi.$ Set $D_\pi:=H(C_A,m+d)$ where $d$ is
the differential in $A$.  The versal Dirac bundle resonates with how
one thinks of the quantum master equation: a solution of the quantum
master equation over a ring $A$ with no differential defines a new
background $B_\pi$ defined over $A$.  By taking the Dirac space of
$B_\pi$ one obtains an $A$ module which is the fiber of a bundle over
the moduli space.

Next, we describe the construction of a flat connection on the versal
Dirace bundle $\widetilde{D}$ of a smooth formal background.
According to Grothendieck, a flat connection on a bundle over a
functor is given by a collection of $A$ module isomorphisms
$D_0\otimes A \to D_\pi$ for each point $\pi$ over $A$, identifying
the fiber over $\pi$ with the fiber over $0$.  It will be convenient
to use a representing couple $(R,\Pi)$ where we take $R$ to have zero
differential.  Using the representing couple, we replace $\widetilde{D}$
with the $R$ module $D_\Pi$ and obtain a flat connection by defining
an isomorphism $D_0\otimes R \to D_\Pi.$ This suffices to define a
flat connection on the versal Dirac bundle since any other point $\pi$
over $A$ is given by a morphism of parameter rings $R\to A$ inducing
$D_\Pi \to D_\pi$ and induces an isomorphism $D_0\otimes A \to D_\pi$.

In terms of a representing couple $(R,\Pi)$, we have the moduli space
$\M=\spec(R)$ and $R=\Lambda[[H^*]]=k[[\hbar,H^*]]$ where $H$ is the
homology of the background and $H^*$ is the linear dual
$H^*=\hom_k(H,k)$.  Since $H=\hom_\Lambda (R,k[t]/t^2)$ is isomorphic
to the derivations of $R$, the vector space $H$ is naturally
isomorphic to the tangent space of the formal space $\M$ at its base
point.  Let $D_0=D(B)$ denote the Dirac space of the original
background, which is the fiber of the Dirac bundle at the base point.
We will prove that $H\otimes_k R \simeq D_\Pi$, which implies that
$D_0\simeq H[[\hbar]]$ and $D_0\otimes R\simeq H \otimes_k R \simeq
D_\Pi$, which defines a flat connection on $\widetilde{D}$.  From another point of
view, bundles over $M$ correspond to $R$ modules.  In particular,
sections of the versal Dirac bundle $\widetilde{D}$ are elements in
the $R$ module $D_\Pi$, which we will also call the versal Dirac
bundle.  The connection that we define, $\nabla^{1}$, defines a map
\begin{align*}
  \nabla^{1}:D_0 \times D_\Pi  &\to D_\Pi \\
  (X,\Y)&\mapsto \nabla^{1}_X \Y
\end{align*}
which is $R$ linear in the first coordinate $\nabla^{1}_{fX}\Y= f
\nabla^{1}_X \Y$ and satisfies the $\hbar$-connection equation in the
second
$$\nabla^{1}_X(f\Y) = \hbar X(f)\Y+(-1)^{|f|} f \nabla^{1}_X
\Y \text{ for $f\in R$, $X\in D$, $\Y\in D_\Pi$}.$$ We prove that
$\nabla^{1}$ is flat:
$$\nabla^{1}_{X}\nabla^{1}_{Y}\mathcal{Z}
-\nabla^{1}_{Y}\nabla^{1}_{X}\mathcal{Z} -
\nabla^{1}_{[X,Y]}\mathcal{Z} =0.$$ Equivalently, we can use
$\Omega:=R[[H^*[1]]]=k[[\hbar, H^*, H^*[1]]]$, the module of
Kahler differentials on $R$, to write the connection using one-forms:
$$\nabla^{1}:D_\Pi \to D_\Pi \otimes_R \Omega^1.$$
In fact, what we define is a ``chain level'' flat connection on the
bundle of cyclic modules $C_R=(P\otimes_\Lambda R)\varphi$
$$\nabla^{chain}:C_R \to C_R\otimes_R \Omega^1$$ which 
descends to $\nabla^{1}$ on the Dirac space.  However, a deeper
analysis of $\nabla^{chain}$ on $C_R$ reveals that when
$\nabla^{chain}$ is compressed from the chain level to the Dirac space, it
manifests as a superconnection $$\nabla:D_\Pi\to D_\Pi \otimes_R
\Omega$$ which can be decomposed as
$\nabla=\nabla^{1}+\nabla^2+\nabla^3+\cdots$ with each $\nabla^i:D_\Pi
\to D_\Pi \otimes \Omega^i.$ In a basis
of $D_\Pi$, $\nabla^i=d_{dR}+A^i$ where $d_{dR}$ is the de Rham
differential in $\Omega$ and $A^i$ is a matrix of $i$-forms.  We prove
that the quantum superconnection $\nabla$ is flat.  The relationship
between $\nabla^{1}$ to $\nabla$ is analogous to a familiar one
concerning differential graded algebras.  Given a dga $(A,d,\cdot)$,
the product descends to 
an associative structure in the homology $H(A,d)$.  However, a
more refined structure which always exists is a minimal
$A_\infty$ structure in $H(A,d)$ that is quasi-isomorphic to the dga
$(A,d,\cdot)$.  The first part of such a minimal $A_\infty$ structure
on $H(A,d)$ coincides the the associative structure in $H(A,d)$, but
the $A_\infty$ structure has more information, it determines the dga
$(A,d,\cdot)$ up to homotopy equivalence.  

The connection $\nabla^1$ is defined canonically from a representing
couple $(R,\Pi)$ but our construction of the superconnection depends
on an $R$ module decomposition $C_R\simeq D'\oplus E \oplus F$  where
$D'\simeq D_\Pi$ are representatives for the Dirac bundle, 
and $m^\Pi:F \to E$ is an isomorphism.  Such a Hodge decomposition of a
chain complex is often used in order to obtain a
minimal $A_\infty$ structure in the homology of a dga.

\subsection{Heisenberg versus Schr\"odinger}
As mentioned above, if $\pi\in \QM_B(A)$ and the differential in $A$
is zero, then we can define a new background with parameters in $A$.
The way it was described above was according to the Heisenberg
representation.
Explicitly,
\begin{equation}
  B_\pi^{Heis}:=(P\otimes A, m^\pi, N\otimes A, \varphi\otimes 1)
\end{equation}
becomes a background: $(m^\Pi)^2=0$ (trivially) and the quantum master
equation says that $m^\Pi\varphi=0$ (nontrivially).  We think of
$B_\pi^{Heis}$ as an evolution of $B$ in which
the structure $m$ that evolves: $$m\leadsto
m^\pi=U^{-1}mU \text{ where }U=e^{-\pi/\hbar}.$$ The structure $m^\pi$ depends on a parameter $t$ in a
parameter ring $R$.  For $t=0$, one has the initial structure
$m^\pi(0)=m$ (since $\pi\in P\otimes \m_R$), and we imagine as $t$
varies, $m^\pi(t)$ varies over the space $\spec(A)$ of structures path
connected to the initial one.  

There is also a Schr\"odinger interpretation of the evolution $B$ to
the background
\begin{equation}
  B_\pi^{Schr}:=(P\otimes A, m, N((\hbar))\otimes A,
  e^{-{\pi}/{\hbar}} \varphi).
\end{equation}
In $\widetilde{B}^{Schr}$, it is the vaccuum $\varphi$, rather than the
structure $m$, that evolves: $$\varphi \leadsto
U\varphi.$$  We permit the module $N$ and the evolving vacuum
$U\varphi=e^{-\pi/\hbar}\varphi$ to
have contain negative powers in $\hbar$, but only finitely many.  
An expansion in terms of 
$A$ coordinates $t_i$
\begin{equation}
  e^{-{\pi}/{\hbar}}\varphi=\sum_{j= 0}^\infty (-1)^j\frac{ \pi^j}
  {\hbar^j j!}\varphi=\sum_{ \substack{ i_1\geq 0,\dots, i_n\geq 0 \\
      \, \\ j\geq -(i_1+\dots+i_n)}} t_{i_1}\dots t_{i_n} \hbar^j p_j^{
    i_1, \dots, i_n} \varphi,
\end{equation}
for $p_j^{i_1,\dots, i_n}\in P$ shows that $e^{-{\pi}/{\hbar}}
\varphi$ is a formal Laurent series for each fixed total degree of the
$t_i$'s. 
Now we prove the equivalence 
of the Heisenberg and Schr\"odinger perspectives
\begin{proposition}The backgrounds $B_\pi^{Heis}$ and $B^{Schr}_\pi$
  are quasi-isomorphic.
\end{proposition}
\begin{proof}
The morphism of quantum
backgrounds $$(\sigma, \tau):B_\pi^{Heis} \to 
B_\pi^{Schr}$$
defined by $$\sigma(\alpha)=U\alpha U^{-1}  
\text{ for }\alpha \in P \otimes A \text{ and
}\tau(\psi)= U \psi \text{ for } \psi\in
N\otimes A.$$
The Heisenberg Dirac module is computed by putting the parameters in
the chain complex and twisting the differential:
$$D(B_\pi^{Heis})=H(C^{Heis}_\pi,m^{\pi})\text{ where }
C^{Heis}_\pi=(P\otimes A)\varphi.$$  In the Schrodinger picture, we twist the
complex and leave the differential unchanged:
$$D(B_\pi^{Schr})=H(C^{Schr}_\pi,m) \text{ where }C_\pi^{Schr}=(P\otimes
A)e^{-\pi/\hbar}\varphi.$$
The map $D(B_\pi^{Heis})\to D(B^{Schr}_\pi)$ induced by
$(\alpha,\tau)$ is an isomorphism  
\begin{gather*}
U^{-1}mUx\varphi=0\Leftrightarrow
mU xU^{-1} U\varphi=0\Leftrightarrow m\sigma(x)\tau(\varphi)=0\intertext{and}
U^{-1}mUx\varphi=y\varphi\Leftrightarrow
mU xU^{-1} U\varphi=U yU^{-1}U \varphi\Leftrightarrow
m\sigma(x)\tau(\varphi)=\sigma(y)\tau(\varphi).
\end{gather*}
\end{proof}

\subsection{The isomorphism $H \otimes_k R \to D_\Pi$}
We now switch to using the Schrodinger representation.  Since we use
this picture for the remainder, we do not use the superscript
$\emph{Schr}$.  Set $U:=e^{-\Pi/\hbar}.$  Then the primary chain
complex consists of the $R$ module
$$C_\Pi=P\otimes R U\varphi$$
with differential given by multiplication by $m$.  
The closed elements correspond to elements $\mathcal{X}\in P\otimes R$
satisfying $m \mathcal{X}U\varphi=0$ and the exact elements
correspond to elements $\mathcal{Y}\in P\otimes R$ satisfying
$\Y U\varphi=m\mathcal{X}U\varphi.$
Since $H$ is the $k$-vector space of $k$-linear derivations of $R$, 
$H$ also acts on $P\otimes R$ as $P$ linear derivations.  
However, a $P$ linear derivation of $R$ maps $C_\Pi\to
\frac{1}{\hbar}C_\Pi$ since differentiating $U$ may introduce one
negative power of $\hbar$.  However $\hbar H$ acts on $C_\Lambda$ by
$P$ linear derivations.

\begin{theorem}\label{connection_thm_1}
  If $B$ is a smooth formal background represented by the couple
  $(R,\Pi)$, then the map $\Phi^0:H\to C_\Pi$ defined by
  $\Phi^0(x)=\hbar x U\varphi$ induces an
  isomorphism $H\otimes_k R \simeq D_\Pi.$
\end{theorem}
\begin{proof}
  Since $\Pi$ is a solution to the master equation, $0=mU\varphi$.
  Apply the derivation $\hbar x$ to get $\hbar x m U \varphi$.  Since
  $m\in P$ and $x$ is a $P$ linear 
derivation of $R$, $m$ and $x$ commute and we
  have $0=m\hbar x
  U \varphi$.  Now, $\hbar x
  (U)\varphi=\hbar x(U) U^{-1}U\varphi$ and since $\hbar x(U)U^{-1}\in
  P\otimes R$ contains no negative powers of $\hbar$, $\Phi^0=\hbar
  x(U)\varphi$ is 
  a closed element in $C_\Pi$.  Taking the homology class of $\Phi^0$
  defines a map $H\to D_\Pi$.

To see that the map, when the coefficients are extended to $R$ becomes
an isomorphism $H\otimes R \to D_\Pi$, 
first observe that it is injective. If $\Phi^0(s)=\hbar x U\varphi$
satisfies $\hbar x U\varphi=m\Y U\varphi$ for some $\Y U\varphi \in
C_\Pi$, then reducing this equation modulo $\hbar$ and $H^*$ implies
that $x$ satisfies $x\varphi = m y\varphi$ for $y=\mathcal{Y} \mod
\hbar, H^*$.

To see that $\Phi^0$ is surjective $D_\Pi$, note that if $m\mathcal{X}
U \varphi=0$ then reducing this equation modulo $H^*$
yields $mX_1\varphi=0$ for $X_1=\mathcal{X} \mod \hbar H^* \in
(P\otimes R)/H^*\simeq P$.  So
$X_1$ defines a class in $H\otimes R/H^*$ and $\Phi^0(X_1)=\mathcal{X} 
\mod H^*.$  Then,
$\mathcal{X}-\Phi^0(X_1) \mod (H^*)^2$ is an $m$ closed
element of $(P\otimes R)/(H^*)^2$ hence defines an element
$X_2\in H\otimes R/(H^*)^2$.  Then 
$\mathcal{X}=\Phi^0(X_1+X_2) \mod (H^*)^3$, and so on.  This
builds inductively an element in $X=X_1+X_2+X_3\in H\otimes
R$ with $\Phi^0(X)=\mathcal{X}\varphi.$
\end{proof}

By reducing the isomorphism $H \otimes_k R \simeq D_\Pi$ modulo $H^*$,
we see that there is no $\hbar$ torsion in $D$.

\begin{corollary}\label{classical-to-quantum}
  If $B$ is a smooth background then $D(B)\simeq
  H(B)[[\hbar]].$  
\end{corollary}
\begin{remark}
In coordinates, say $\{x_i\}$ is a basis for $H$ with dual basis 
$\{t_i\}$, we have $x_i$ acts by $\hbar \frac{\partial}{\partial t_i}$
on $R=k[[\hbar,t_i]]$ and $P[[\hbar, t_i]].$   A versal
solution to the quantum master equation $\Pi$ has the form
\begin{multline*}
\Pi=\sum  \left( x_i t_i+ \hbar x^{(1)}_i t_i + \hbar^2 x^{(2)}_i t_i
  + O(\hbar^3) \right) \\
+  \left(x_{ij} t_it_j + \hbar x_{ij}^{(1)}
  t_i t_j \hbar^2 x^{(2)}_i t_j+O(\hbar^3)\right)+\cdots
\end{multline*}
with $x_{i_1, \ldots, i_n}\in P$.  The first terms
$\{x_i\}$ are the $k$-basis for $H$ and the linear in $t$ terms 
$\{X_i:=x_i+ \hbar x^{(1)}_i + \hbar^2 x^{(2)}_i 
  + \hbar^3\cdots\}$ are a $k[[\hbar]]$ basis for $D$ and
$\left\{\mathcal{X}_i:=\hbar \frac{\partial e^{-\Pi/\hbar}}{\partial
    t_i}\right\}$ is an 
$R=k[[\hbar, H^*]]$-basis for $D_\Pi$.  In the case that $\Pi$
commutes with its derivatives in $P$, $\mathcal{X}_i=\frac{\partial
  \Pi}{\partial t_i}.$  The isomorphisms $H\otimes \Lambda \to D$ and
$H\otimes R \to D_\Pi$ are given by $x_i \mapsto X_i$ and $x_i
\mapsto \mathcal{X}_i.$
\end{remark}

\begin{remark}
Recent
developments \cite{js3,jt2} 
suggest that the converse to 
Corollary \ref{classical-to-quantum} is true and the concept of
special coordinates in string theory is related to the way in which an
isomorphism $D(B)\simeq H(B)[[\hbar]]$ can be extended to produce a
``special'' versal solution to the quantum master equation.
\end{remark}

So, as described in the conecptual description of the connection,
Theorem \ref{connection_thm_1} implies that $D_0\otimes R\simeq
H\otimes_k R \simeq D_\Pi$ and therefore defines a flat connection on
$D_\Pi.$ Now we translate from the Grothendieck formal geometry
picture to the covariant derivative picture of a flat connection.  The
translation amounts to observing that $D_0=H[[\hbar]]$ 
is the $\Lambda=k[[\hbar]]$ module of $\Lambda$ linear derivations of 
$R=k[[\hbar,H^*]]$, hence act as $P$ linear derivations of $C_\Lambda$
and $D_\Lambda$.
\begin{definition}\label{chain-connection}
  Define $\nabla^{chain}: D_0 \times C_\Pi \to C_\Pi$
by $(X,\Y)\mapsto \nabla^{chain}_X(\Y)=\hbar X(\Y)U^{-1}.$
\end{definition}
First note that $\hbar X (\mathcal{Y}U\varphi)=\hbar
X(\Y U)U^{-1}U\varphi$ and that $\hbar X(\Y U)U^{-1}$ has no negative
powers of $\hbar$.  This proves that $\nabla_X:C_\Pi \to C_\Pi$.
Now, since $m\in P$ and $X$ acts as a $P$ linear derivation, 
$\nabla_X(m\Y)=X(m\Y)=mX\Y$ so $\nabla_X$ is a chain map, hence
defines a map
\begin{equation}
  \nabla^{1}: D_0 \times D_\Pi \to D_\Pi
\end{equation}

\subsection{Superconnection on the versal Dirac
  module}
First, refashion $\nabla^{chain}$ in terms of differential forms.
Let
$\Omega^k=R[[H^*]]\otimes_R S^k(H^*[1])$ be the module of $k$
forms on $\spec{R}$ and $\Omega:=\Pi_k\Omega^k$
be the module of Kahler differential forms.  
Let $d_{dR}:\Omega \to \Omega$ denote the de Rham differential; i.e.,
the shift functor $H\mapsto H[1]$ extended
to a $\Lambda$-linear derivation of $\Omega$.  Therefore
$\hbar d_{dR}:C_\Pi\otimes_R \Omega \to C_\Pi\otimes_R \Omega $ is a square
zero $P$ linear derivation.  Since $m\in P$, $[m,\hbar d_{dR}]=0$. 
Therefore, $d_{dR}$ descends to give a map
$$\nabla^1:D_\Pi \to D_\Pi \otimes_R \Omega^1$$
This is precisely the same as the connection defined before, expressed
in terms of a connection one forms.

However, one can do better.  We
now explain how to transfer $d_{dR}$ to a superconnection 
$\nabla: D_\Pi \to D_\Pi\otimes_R \Omega.$  In order to do so,
we use one more piece of data.
Consider $\Phi^0:H \to C_\Pi$ as in Theorem
\ref{connection_thm_1} defined by $\Phi^0(x)=\hbar x(U).$  
Since 
$\Phi^0(H)$ is closed in $C_\Pi$ and $R\Phi^0(H)\simeq D_\Pi$, we
may choose an acyclic complement $E$ in $C_\Pi$ so that 
\begin{equation}\label{directsum}
C_\Pi=\Phi^0(H)\oplus E.
\end{equation}
Now, the equation 
\begin{equation}\label{eq1}
m\Phi^0=0
\end{equation}
 implies $m \hbar d_{dR} \Phi^0=0$.  Therefore,
\begin{equation}\label{eq2}
\hbar d_{dR} \Phi^0 = A^1 \Phi^0 + m \Phi^1
\end{equation}
where $A^1:D_\Pi \to
D_\Pi \Omega^1$ and $\Phi^1:H \to E\otimes \Omega^1$.  
Since $d_{dR}^2=0$, applying $\hbar d_{dR}$ to Equation \eqref{eq2} yields
$$(\hbar d_{dR}A^1) \Phi^0 - A^1 \hbar d_{dR}\Phi^0 - m \hbar d_{dR}\Phi^1=0.
$$
Using Equation \eqref{eq2} again get
\begin{multline*}
0=(\hbar d_{dR}A^1) \Phi^0 - A^1 (A^1 \Phi^0 + m \Phi^1) -m
\hbar d_{dR}\Phi^1\\ =(\hbar d_{dR}A^1-A^1A^1 )\Phi^0 + m(\hbar d_{dR}\Phi^1-A^1 \Phi^1).
\end{multline*}
This one equation splits into the two equations
\begin{gather}
\hbar d_{dR}A^1-A^1A^1=0\\
m(\hbar d_{dR}\Phi^1-A^1 \Phi^1)=0
\end{gather}
Now, the second equation defines a cycle in $C_\Pi\otimes_R \Omega$ and
so implies that $\hbar d_{dR}\Phi^1-A^1 \Phi^1=A^2
\Phi^0 + m \Phi^2$ for unique $A^2:D_\Pi \to D_\Pi \otimes_R \Omega^2$
and $\Phi^2:H \to E\otimes_R \Omega^2.$
Iterating this procedure proves the following theorem:
\begin{theorem}\label{conn-A}
  For every versal solution $\Pi$ to the quantum master equation, and
  every choice of acyclic complement $E$ of $\Phi^0(H)$ in $C_\Pi$, 
  there exist unique maps $A^n:D_\Pi \to D_\Pi\otimes \Omega^n$
  for $n\geq 1$ and $\Phi^n:H\to E\otimes \Omega^n$ for
  $n\geq 0$ satisfying
  \begin{align}
        \label{d(A)}
    \hbar d_{dR}A^n&=\sum_{j=1}^n (-1)^{n+j} A^{n+1-j}A^{j}, \\
    \label{d(phi)}
    \hbar d_{dR}\Phi^n_i&=\sum_{j=0}^n A^{n+1-j}\Phi^j+m\Phi^{n+1}.
  \end{align}
\end{theorem}
\begin{definition}\label{conn_def}
  Define $\nabla:D_\Pi \to D_\Pi\otimes_R \Omega$ by 
  $$\nabla:= \hbar d_{dR}+A^1+A^2+\cdots$$ 
\end{definition}
\begin{theorem}\label{conn_thm}
  $\nabla$ is flat.
\end{theorem}
\begin{proof}
  Notice, that $\nabla$ satisfies $$ \nabla(f\Y)=\hbar d_{DR}(f)\wedge
  \Y+(-1)^{|f|} f\nabla(\Y) $$ for any $f\in R$ and $\Y\in D_\Pi$.
  Equation \eqref{d(A)} implies $\nabla^2=0$, which
  is the flatness condition.
\end{proof}

\begin{remark}
  As discussed in the next 
  section, this small connection carries all of the information
  about the correlation functions (except for a choice of one-point
  functions).  However, $\nabla$ carries much more information, what
  might be thought of as homotopy correlation functions.  For example,
  in the $B$-model, the homotopy correlation functions should give a
  chain level generalization of the Frobenius manifold structure
  already discovered \cite{B2, BK}.
\end{remark}

\section{Path integral, correlation functions, and generalizations}
In this section, we relate the quantum superconnection to the
correlation functions of the quantum field theory defined by the
quantum background $B$.  We assume $B$ is finite, but not necessarily
smooth formal.

\subsection{Path integral}

\begin{definition}\label{path-integral}
  Let $B=(P,m,N,\varphi)$ be a background.  We define a (chain level)
  path integral pairing for $B$ to be a $k[[\hbar]]$ module map
  $\int:N \to k[[\hbar]]$, which we denote by $\psi \mapsto \int
  \psi,$ satisfying
  \begin{itemize}
  \item ($P$ Module axiom) $\ds \int (\alpha\beta)\psi = \int \alpha
    (\beta \psi)$ for all $\alpha, \beta \in P$, and $\psi \in N$,
  \item (Stokes axiom) $\ds \int m\psi =0$ for every $\psi \in N$.
  \end{itemize}
\end{definition}
We call the condition $\ds \int m\psi =0$ Stokes axiom because, if we
use $\int_\psi \alpha$ to denote $\int \alpha \psi$, the condition
$\int_\psi m=0$ implies
$$\int_{\partial \psi} \alpha = \int_\psi D\alpha , $$
where $\partial$ and $D$ are differentials in $P$ and $N$ defined by
$$\partial(\psi)=m\psi \text{ and }D\alpha = m \alpha.$$

Given a finite background $B$ and a ring $(R,\delta)$ which represents
$\QFT_B$, one can extend the module $N$ to $C_\Lambda$ which
is a left $P$ module with the action of $m$ extended to be that of
$m\otimes 1.$
A path integral pairing $\int$ for $B$ extends
$$C_\Lambda  \to R((\hbar))$$ and
extends with module compatibility property and the Stokes property.

\begin{definition}
  If $\Pi$ is a versal solution to the quantum master equation, then
  $e^{-{\Pi}/{\hbar}}\varphi\in C_\Lambda \otimes R$ and we define
  the partition function $\mathcal{Z}\in R((\hbar))$ by
  $$\mathcal{Z}=\int e^{-{\Pi}/{\hbar}}\varphi.$$
\end{definition}

\subsection{Generalized Ward identities}
For a couple $((R,\delta),\Pi)$ representing $\QFT_B$, we have the
master equation
$$me^{-{\Pi}/{\hbar}}\varphi -\hbar\delta(e^{-{\Pi}/{ \hbar}} )\varphi=0.$$
Integrating (and Stokes axiom) gives
$$\int \hbar \delta (e^{-{\Pi}/{ \hbar}})\varphi=\hbar\delta \int e^{-{\Pi}/{\hbar}}\varphi=0.$$
The integral is linear over $R$, giving the following differential
equation for the partitition function:
\begin{equation}\label{wardequation}
  \hbar\delta \mathcal{Z}=0.\end{equation}
This non-trivial equation captures the most general 
symmetries of the integral, including, in some example, the Ward
identities.  In the case that the background is smooth, the Ward
identities vanish.  The conclusion is
\begin{quote}
\emph{The obstructions to deforming the background
manifest themselves as symmetries of the partition function}.
\end{quote}

\subsection{Correlation functions}
Now suppose that $B$ is smooth formal.  Choose $R$ so that $\delta=0$.
The tangent space $H=H(B)$ plays the role of the classical observables and
the Dirac space $D=D(B)$ plays the role of physical observables.   In
the smooth case, there is an isomorphism $H[[\hbar]]\to D$, hence we
may think of the two classes of observables interchangeably.  We
think of the Dirac bundle $D_\Pi$ as consisting of the observables for all
theories in a neighborhood of $B$ in the moduli space.  
\begin{definition}
  For each $n=0,1,2,\ldots $, we define multilinear maps, called
  $n$-point correlation functions,
$$\<\cdots  \>:H^{\otimes n}\to R((\hbar))$$ by 
$$Z_{i_1}\otimes \cdots \otimes Z_{i_n}\mapsto\<Z_{i_1}, \ldots,
Z_{i_n}\>:= 
 \int \hbar^nZ_{i_n}\cdots Z_{i_1} e^{-{\Pi}/{\hbar}}\varphi.$$
\end{definition}
Here, the technique
of computing correlation functions by differentiating a family of
action functionals here becomes a definition since each observable
$Z_i\in H$ is a derivation of the representing ring $R$.
The $n$ point correlation functions are completely determined by the
one point correlation functions and the small connection
$\nabla^{1}$; i.e., the one-form part of $\nabla$.  To see this note
that $\<Z_i\>=\<\Phi^0(Z_i)\>$.  Applying $\hbar d_{dR}$ gives a one-form
$\hbar d_{dR}(Z_i)$ which when evaluated on the tangent vector $Z_j$ gives
$$\< \hbar d_{dR}(\Phi^0(Z_i))(Z_j)\>=\<Z_j,Z_i\>.$$
The equations in Theorem 
\ref{conn-A} imply 
$$ \<Z_j,Z_i\>=\<A^1\Phi^0(Z_i)+m\Phi^1(Z_i)\>=A^1\<\Phi^0(Z_i)\>.$$
Next, one finds that 
$$\<Z_k,Z_j,Z_i\>=\hbar d_{dR}A^1(Z_k)\<\Phi^0(Z_i)\>+A^1A^1\<\Phi^0(Z_i)\>.$$
By recursively differentiating and eliminating the exact terms, one
can compute all correlation functions in terms of the $A^i$, and their
derivatives.

\subsection{Homotopy correlation functions}
A chain level path integral pairing gives rise to a linear functional
on the versal Dirac bundle and this linear functional on
$D_\Pi$ is sufficient to determine the correlation
functions.  So, one may propose to define a path integral, or a
``cohomological path integral'' to contrast it with a chain level path
integral, as any linear functional $F\in \hom(D_\Pi,R)$.
This avenue will lead to the definition of the one-point correlation
functions as $\<[Z_i]\>=F([Z_i])$ and the definition of the $n$-point
correlation functions by equations already given.
However, there is information, beyond
correlation functions, that a chain level path integral pairing can
detect that is invisible to any cohomological path integral.

To explain, we first make the correlation function into
a function of the moduli by $\int \Phi^0$.  
Then, the $n$
point correlation functions are obtained as derivatives of the one
point correlation functions and are inductively determined by
integrating Equation \eqref{d(phi)} for $n=0$:
\begin{equation}\label{one-form-recursion}
  \hbar d_{dR} \int \Phi^0 = A^1\int\Phi^0.
\end{equation}
Thus, one might as well summarize all of the $n$ point correlation
functions into the single fundamental correlation function
\begin{equation}\label{correlation-form}
  \int \Phi^0 
\end{equation}
with the understanding that the correlation function satisfies
Equation \eqref{one-form-recursion}.

Note that the information contained in the boundary  term $m\Phi^1$
from equation \eqref{d(phi)} is lost after integration, but may be
retained as a correlation one-form
$$\int \Phi^1.$$ 
Think of the correlation one-form as a homotopy $1$-point correlation
function, from which many homotopy $n$-point functions can be derived
by differentiation.  These homotopy correlation functions that may be
derived from $\int \Phi^1$ are summarized inductively by integrating
equation \eqref{d(phi)}:
$$ \hbar d_{dR}\int \Phi^1=A^1 \int \Phi^1+A^2 \int \Phi^0.$$
This discussion suggests an efficient way to handle these correlation
and homotopy correlations by defining the correlation $k$ form by
$\int \Phi^k$ which satisfies $d_{dR}\int \Phi^k=A^1 \int \Phi^{k}+A^2
\int \Phi^{k-1}+\cdots+A^k\int \Phi^0.$ Better yet:
\begin{definition}We define the primary correlation form to be the
  function  $$\int \Phi:H \to \Omega $$ where
  $\Phi=\Phi^0+\Phi^1+\Phi^2+\cdots. \in C_\Pi \otimes_R \Omega$ 
\end{definition}
\begin{theorem}The primary correlation form satisfies $\ds \nabla \int
  \Phi =0.$
\end{theorem}

\section{Application---dBV algebras}\label{Application}
There are many situations which give rise to quantum backgrounds.  In
each, the resulting quantum flat superconnection provides a
penetrating tool for investigating the situation.  Here, we give one
example which we feel may elucidate the many definitions in this
paper.  We begin with a differential BV algebra and construct a
quantum background from it.  Then, we indicate how the various
features of the background correspond to the dBV algebra.  In
particular, the quantum superconnection derived from the $\QFT$
functor deeply probes the homotopy theory of the BV algebra.  Related
ideas appear in \cite{js1}.

Let $(V,d,\Delta,\wedge)$ be a differential BV algebra.  This means
that $V$ is a graded vector space and $d, \Delta$ are commuting differentials on
$V$ with the properties
\begin{itemize}
\item $(V,d,\wedge)$ is a differential graded commutative,
  associative, algebra, and
\item $\left(V[-1],d, (\,,\,)\right)$ is a differential graded Lie
  algebra, where the bracket is defined by
  $(v,w):=(-1)^{|v|}\Delta(v\wedge w)-(-1)^{|v|}\Delta(v)\wedge
  w-v\wedge \Delta(w)$ for homogeneous vectors $v$ and $w$.
\end{itemize}
We use parentheses $(\,,\,)$ for the Lie bracket, reserving square
brackets $[\,,\,]$ always for the graded commutator.  For convenience,
we assume that the dBV algebra is fairly simple; assume $V=SU$ for a
finite dimensional graded vector space $U$ for which $\wedge$ is the
associative symmetric product in $SU$.  This assumption allows us to
operate easily in coordinates.  Let $\{q_1, \ldots, q_n\}$ be a
homogeneous basis for $U$.  Then elements of $V$ are polynomials in
the variables $\{q_i\}$, the wedge product is the ordinary graded
commutative product of polynomials, and we may abbreviate $q_i\wedge
q_j$ by $q_iq_j$.  The operator $d$, as a derivation of the product,
is a first order differential operator, and $\Delta$, as a BV
operator, is a second order differential operator.  Such operators
have expressions in the basis $\{q_i\}$ as
\begin{gather}\label{def_of_d}
  d=\sum_{i=1}^n a_i(q) \frac{\partial}{\partial q_i} \intertext{ and
  } \Delta=\sum_{i,j=1}^n b_{ij}(q)\frac{\partial^2}{\partial q_j
    \partial q_i}\label{def_of_delta}
\end{gather}
$a_i(q)\in V$ for each $i$ and $b_{ij}(q)\in V$ with
$b_{ij}(q)=-(-1)^{|q_i||q_j|}b_{ji}(q)$ for each $i$ and $j$.  From
such a dBV algebra, we now define a background
$B_{V,d,\Delta,\wedge}=(P,m,N,\varphi)$.

\subsection{The ring $P$.}  We define the ring $P$ as $P=W(U)$, a
graded Weyl algebra on the vector space $U$, defined to be
\begin{gather*}
  W(U):=T(U\oplus U^*)[[\hbar]] \big\slash J \intertext{ where $J$ is
    the left ideal generated by} [q,q'],\quad [p,p'],\text{ and }
  [q,p]=\hbar p(q).
\end{gather*}
for any $q, q' \in U$, and $p, p' \in U^*$, and $[\,,\,]$ is the
graded commutator.  The product in $P$ is induced by the tensor
algebra.  As a $k[[\hbar]]$ module, $P\simeq (S(U\oplus
U^*))[[\hbar]]$ and as a $k$-algebra $P/\hbar P \simeq S(U\oplus
U^*).$ In coordinates, say $\{q_i\}$ is a basis for $U$ and $\{p^i\}$
is the dual basis of $U^*$, each elments of $P$, which is an
equivalence class in $T(U\oplus U^*)[[\hbar]]$ is represented uniquely
by a polynomial in the $\{q_i\}$ and $\{p^i\}$ in normal ordering with
``all the $p$'s on the right.''  Multiplication is carried
mechanically out by concatinating the two polynomials and using the
commutation relation $p^iq_j=(-1)^{|q_j||p^i|}q_jp^i-\hbar\delta^i_j$
repeatedly until obtaining a polynomial with the $p$'s on the right,
thereby obtaining the unique representative for the product in $P$.

The Weyl algebra $P=W(U)$ may be familiar as the ``canonical
quantization'' of the symplectic vector space $U\oplus U^*$.  The ring
$S(U\oplus U^*)$ can be thought of as the ring of (polynomial)
functions on a symplectic space, with the Poisson bracket defined by
$\{q_i,p^j\}:=\delta_i^j$.  Then $P$, isomorphic as a $k[[\hbar]]$
module to $S(U\oplus U^*)[[\hbar]]$, is a one-parameter deformation of
$S(U\oplus U^*)$ over $k[[\hbar]]$ with the infinitesimal deformation
being the Poisson bracket.  The product we defined in $P$ might be
called a ``star product.''

While it is easy (and tedious) to perform calculations in $P$ using
explicit polynomials in $p$ and $q$, one can understand the elements
in $P$ and the product in $P$ in a coordinate free way.  Observe that
$$P/\hbar P\simeq S(U\oplus U^*) \subseteq  V\otimes V^* = \hom_k(V,V),$$
and hence, as a $\Lambda$ module,
\begin{equation}\label{weyl_hom}
  W(U)\subseteq
  \hom_k(V,V)[[\hbar]].\end{equation}
So, without choosing a basis for $U$, one can regard
elements of $P$ as power series in $\hbar$ with
coefficients that are the finite $k$-linear operators on $V$;
i.e. finite sums of operators $S^jU \to S^kU$. 
The
identification in equation \eqref{weyl_hom} is as free
$\Lambda=k[[\hbar]]$
modules;
to acquaint oneself with the multiplication in $P$ from this point of
view, we expand
$$\alpha\beta = \alpha\wedge \beta + \hbar
\{\alpha,\beta\}+\mathcal O(\hbar^2)$$ for $\alpha,\beta\in
\hom_k(V,V)$.  We interpet some of the $\hbar$ terms on the right for
a couple of illustrative cases:
\begin{example}
  Suppose $\alpha:S^jU\to S^kU$ and let $\beta:S^iU \to S^jU$. Then,
  \begin{equation}
    \alpha\beta = \alpha\wedge \beta+(\hbar \text{ terms})+(\hbar^2 \text{
      terms})+\cdots +\hbar^j (\alpha \circ \beta)
  \end{equation}
  where $\alpha\wedge \beta: S^{i+j}U \to S^{j+k}U$ is the graded
  commutative wedge product of linear tranformations and $\alpha \circ
  \beta: S^i U\to S^kU$ is composition in $\hom(V,V)$.  Hence, $\star$
  extends the commutative associative product on $\hom(V,V)$ toward
  the noncommutative composition of homomorphisms.  The $\hbar^r$
  terms for $r=2,3, \ldots,j-1$ can be understood in terms of some
  multiple-gluing operations using a topological model.  In this model
  $\hbar$ is related to genus or Euler characteristic \cite{DTT}.
\end{example}
\begin{example}
  Let $\alpha:S^iU\to U$ and let $\beta:S^jU \to U$.  Then,
  \begin{equation}\label{item2} [\alpha,\beta]= \hbar
    \underline{\left[\overline{\alpha}, \overline{\beta}\right]_G}
  \end{equation}
  where the bracket on the left is the graded commutator in $P$ and
  the bracket $[\,,\,]_G$ on the right is understood by way of the
  Gerstenhaber bracket in the $\coder(V)$, the Hochschild complex of
  the algebra $(V[-1], d, \wedge)$.  That is, to compute
  $[\alpha,\beta]$ in $P$, one lifts $\alpha$ and $\beta$ to
  coderivations $\overline{ \alpha}, \overline{\beta}: V\to V$, takes
  the graded commutator of those two coderivations to obtain a
  coderivation $\left[ \overline{ \alpha}, \overline{\beta}\right]_G:
  V \to V$, with the final result of $[\alpha,\beta]$ being $\hbar$
  times the component $S^{i+j-1}U\to U$ that determines the
  coderivation $\left[\overline{\alpha}, \overline{\beta}\right]_G$.
\end{example}

\subsection{The structure $m$.}  We define $m\in P^1$
$$m=\sum_{i=1}^n
a_i(q)p^i +\sum_{i,j=1}^n b_{ij}(q)p^ip^j$$ where $a_i$ and $b_{ij}$
are as in equations \eqref{def_of_d} and \eqref{def_of_delta}.  The
fact that $(V,d,\Delta,\wedge)$ is a BV algebra implies that as an
element of $P$, $m^2=0$.

\subsection{The module $N$.}  We have $U^*\subset P$ and the left
ideal $PU^*$.  The quotient $P/PU^*$ is naturally a $P$ module and we
define $N=P/PU^*$.  Use $|\alpha\>$ to denote the image of $\alpha \in
P$ in $P/PU^*$.  There is an isomorphism of $k[[\hbar]]$ modules
$V[[\hbar]]\simeq N$.  In coordinates, this isomorphism is simple: an
element of $V[[\hbar]]$ is a polynomial in $q$ which represents a
class of polynomials in $P$, and also a class in $P/PU^*$ since $V$
has no $p$'s.  To illustrate, a typical element of $P$ might be
$\alpha=2\hbar q_1q_2+q_1p^2p^3+\hbar^2 p^3$ and typical element of
$N$ might be $\gamma=|2\hbar q_1-q_1q_3\>$.  Then (say $q_1$ and $q_3$
are even and $q_2$ is odd),
\begin{align*}
  \alpha \cdot \gamma &=|(2\hbar q_1q_2+q_1p^2p^3+\hbar^2 p^3) (2\hbar
  q_1-q_1q_3)\>\\
  &=|(4\hbar^2 (q_1)^2q_2-2\hbar (q_1)^2q_2q_3+2\hbar (q_1)^2 p^2
  p^3-\hbar (q_1)^2p^2+\hbar^2q_1q_3p^3-\hbar^3q_1\>\\
  &=|(4\hbar^2 (q_1)^2q_2-2\hbar (q_1)^2q_2q_3-\hbar^3q_1\>.
\end{align*}

\subsection{The vacuum $\varphi$}
The vaccum $\varphi$ is defined to be $|1\>$.  Note that
$m\varphi=m|1\>=|m\>=|0\>$, as required.

\subsection{Relating the background and the dBV algebra}
We begin with a summary of the results.  It will be convenient to
state the summary in terms of several algebras built from
$(V,d,\Delta, \wedge)$:
\begin{align*}
  C&:=(V,d,\wedge) \text{ is a commutative dga over $k$}\\
  L&:=(V[-1],d,(\,,\,) \text{ is a dgLa over $k$}\\
  L^\hbar&:=(V[-1][[\hbar]], d+\hbar\Delta, (\,,\,)) \text{ is a dgLa
    over
    $k[[\hbar]]$}\\
  H&:=\ker(d:V \to V)/\im(d:V\to V) \text{ is
    a vector space over $k$}\\
  H^\hbar&:=\ker(d-\hbar\Delta:V[[\hbar]] \to
  V[[\hbar]])/\im(d-\hbar\Delta :V[[\hbar]]\to V[[\hbar]]) \text{ is a
    $k[[\hbar]]$ module.}
\end{align*}
Note that $H\simeq H(C)\simeq H(L)[1]$.  By the minimal model theory
for $L_\infty$ algebras, there is an $L_\infty$ structure on $H[-1]$,
unique up to quasi-isomorphism, that is quasi-isormorphic to $L$.
Note also, $H^\hbar\simeq H(L^\hbar)[1]$, so there is an $L_\infty$
structure (over $k[[\hbar]]$) on $H^\hbar[-1]$ quasi-isomorphic to
$L^\hbar$.

On the quantum background side, we denote $B_{V,d,\Delta,\wedge}$
simply by $B$, and we have
\begin{align*}
  (R&=k[[H^*,\hbar]],\delta):=\text{ a ring representing $\QFT_B$}\\
  \Pi&\in \QM_B(R,\delta):=\text{ versal solution to quantum master
    equation}\\
  (S&=k[[H^*]],\delta_0):=(R\otimes_\Lambda k,
  \delta\otimes_\Lambda k) \text{
    is the ``classical'' $\hbar=0$ part of $R,\delta$}\\
  \Pi_0&\in \QM_B(S,\delta_0):=\text{ the image of $\Pi$ in $S$ is
    the ``classical'' part of $\Pi$}.
\end{align*}
The most immediate relationships between the background $B$ and the
dBV algebra $(V,d,\Delta, \wedge)$ are summarized in Figure 1.

\begin{figure}[h]
  \label{table1}
  \begin{center}
    \begin{tabular}{|l|l|}
      \hline
      $H(B)=$ the tangent space to $\QFT_B$ & $H$\\
      \hline
      $D(B)=$ the Dirac space of $B$ & $H^\hbar$\\
      \hline
      The ring $S$ &  $k[[ H^*]]$ \\
      \hline
      The ring $R$ &  $k[[ \hbar, H^* ]]$ \\
      \hline
      $(S,\delta_0)$ &  \parbox{6cm}{\smallskip
        the dual of an $L_\infty$ minimal model of $L$\smallskip}
      \\
      \hline
      $\Pi_0$ & \parbox{6cm}{\smallskip a minimal model map $H\to L$\smallskip}
      \\
      \hline
      $(R,\delta)$ &  \parbox{6cm}{\smallskip
        the dual of an $L_\infty$ minimal model of 
        $L^\hbar$\smallskip}\\
      \hline
      $\Pi$ & \parbox{6cm}{\smallskip a minimal model map $H^\hbar \to L^\hbar$\smallskip}
      \\
      \hline
    \end{tabular}
    \caption{Comparing $B$ to $(V,d,\Delta, \wedge)$}
  \end{center}
\end{figure}

In the case that $B$ is smooth formal, we also have $L$ and $L^\hbar$
smooth formal and the moduli space for $\QFT_B$ and $\Def_L$ are the
identified.  On the BV side, we may define the following algebras,
which can be thought of as families of algebras fibered over the
moduli space.
\begin{align*}
  C_t&:=(V\otimes S,d+(\Pi_0,\,),\wedge) \text{ is a commutative dga
    over $k$}\\
  L_t&:=(V[-1]\otimes S,d+(\Pi_0,\,),(\,,\,) \text{ is a dgLa over $k$}\\
  L_t^\hbar&:=(V[-1][[\hbar]]\otimes R, d+\hbar\Delta +(\Pi,\,),
  (\,,\,)
  \text{ is a dgLa over $k[[\hbar]]$}\\
  H_t&=\ker(d+(\Pi_0,\,):V\otimes S\to V\otimes
  S)/\im(d+(\Pi_0,\,))
  \text{ is an $S$ module}\\
  H^\hbar_t&=\ker(d-\hbar\Delta+(\Pi,\,):V[[\hbar]]\otimes R\to
  V[[\hbar]]\otimes R)/\im(d+(\Pi_0,\,)) \text{ is an $R$ module}
\end{align*}
We note that $H_t\simeq H(C_t)$.  So, $H_t$ is a commutative
associative algebra.  Furthermore, by the minimal model theory, $H_t$
has a $C_\infty$ structure, quasi-isomorphic to $C_t$.

On the quantum background side, we have the quantum superconnection
\begin{align*}
  \nabla&=d_{dR}+A^1+A^2+\cdots:=\text{ the quantum superconnection}\\
  \nabla^{1}&=d_{dR}+A^1:=\text{ the small quantum connection}\\
  \nabla_0&=d_{dR}+A_0^1+A_0^2+\cdots:=\text{ the ``classical''
    $\hbar=0$
    part of $\nabla$ }\\
  \nabla_0^{1}&=d_{dR}+A_0^1:=\text{ the ``classical'' $\hbar=0$  part of $\nabla^{1}$}\\
\end{align*}
We conjecture that there are special coordinates for any smooth formal
background.   What that means in the present example of a dBV algebra
is that there is a particular solution $\Pi$ to the quantum master
equation for which the associated flat quantum connection $\nabla$ has
no $\hbar$ dependence.  This is proved in \cite{js3} in a certain
semiclassical case.  When $\nabla^{1}$ does not depend on $\hbar$,
the flatness equation expressed in terms of the conncection one form
$A$
$$\hbar d_{dR}A+A^2=0$$ decouples into two 
equations $d_{dR}A=0$ and $A^2=0$.  The equation $A^2=0$ together with the
torsion-freeness of $\nabla^{1}$ implies that $A$ defines a family
of commutative, associative algebra structures on $H$ parametrized by
$H$.  The condition that $dA=0$ gives additional constraints on the
way theses algebra structures vary with their parameters, implying
that $H$ has the structure of a Frobenius manifold
without a metric \cite{Manin}.  We cojecture that the
superconnection in special coordinates affords $H$ with the structure
of a minimal algebra over an appropriate resolution of the Frobenius manifold
structure.  

\subsection{More general Weyl-type backgrounds}
We illustrate so many details about a background arising from a dBV
algebra because the situation may be familiar, but we emphasize that
it is only an example.  Notice that $m\in P$ obtained from $d$ and
$\Delta$ in a dBV algebra is quite special: it is quadratic in $p$ and
has no dependence on $\hbar$.  One may more generally consider a
background $B=(P,m,N,\varphi)$ where $P=W(U)$ is the Weyl algebra on a
graded vector space $U$, $m$ is any element, at least linear in $p$
(in order to annihilate the vacuum), satisfying $m^2=0$,
$N=SU[[\hbar]]$, and $\varphi=|1\>$.  There seem to be many
interesting examples \cite{DTT}.  Furthermore, for a fixed $m$, one may consider
different modules, possibly highlighting different aspects of the same
structure.

\section{Prospectus}
In conclusion, let us make some very brief remarks about future
directions, to which we are now turning.

In this paper, we constructed a flat quantum superconnection over the
moduli space if the background is smooth formal.  Our attention is now
focused on the non smooth formal case.  In the general case, we intend
to construct a quantum connection $\nabla$ which interacts with the
differential $\delta$ from proposition \ref{main_thm}.  This seems to
be necessary to develop a good minimal model theory for quantum
backgrounds.

We imagine the ideas in this paper can be applied in two different
ways.  The first application is rather direct: apply the framework
described here to mathematical situations that fit.  One noteworthy
example arises in symplectic field theory.  One might summarize the
output of symplectic field theory as an element of square zero in a
particular noncommutative ring.  In other words, the output of
symplectic field theory is a rather good match to the input data of a
background.  It would be interesting to construct the $\QFT$
moduli space and quantum connection $\nabla$ in this example and
interpret these structures in terms of symplectic topology.

A second application is a kind of quantized deformation theory.  By
this we mean start with a classical mathematical structure and produce
a quantum background $B$ whose classical limit is the $L_\infty$
algebra $L$ controlling the given structure's classical deformations.
By ``classical limit'' we mean that $B$ and $L$ share a relationship
much the same as the relationship between $B$ and $L$ exemplified in
Section \ref{Application}.  Then, the $\QFT_B$ moduli space enriches
the classical moduli space in the $\hbar$ direction.  One advantage is
that the superconnection $\nabla$ when expressed in special
coordinates, which are invisible from the classical
deformation theory, is expected to encode invariants of the
original structure.

\end{document}